# NeSTR: Neural $S$-Transform Reconstruction for Gaussian SPDEs

Nacira Agram[1] Fred Espen Benth[2] Jan Rems[3]

September 19, 2026

**Abstract**

We introduce Neural $S$-Transform Reconstruction (NeSTR) for Gaussian stochastic partial differential equations (SPDEs). The learned object is not the random solution itself, but its $S$-Transform on a finite dimensional subspace of the white noise test space. For Wick type SPDEs, NeSTR yields a deterministic parametric PDE for the restricted $S$-Transform, and the stochastic solution is recovered from the Taylor coefficients of this transform at the origin. Thus the neural approximation is tied directly to the Wiener chaos expansion of the solution. The main novelty of NeSTR is an inverse $S$-Transform reconstruction principle combined with finite mode learning: we prove consistency of the recovered chaos coefficients, separate noise mode, chaos truncation, deterministic solver, and neural approximation errors, and give a deterministic learning formulation for the restricted transform. Numerical tests on an additive stochastic heat equation and a multiplicative Wick heat equation show that polynomial spectral representations in the transform variables provide stable access to higher order chaos coefficients. The resulting reconstruction separates the equation dependent deterministic coefficient fields from a reusable Gaussian–Hermite basis. Consequently, NeSTR learns an analytic generating function whose derivatives recover the Wiener chaos coefficients, rather than learning stochastic trajectories directly.

## Contents

[1]Department of Mathematics, KTH Royal Institute of Technology and Digital Futures, 100 44 Stockholm, Sweden. Email: nacira@kth.se. Supported by the Swedish Research Council grant (2020-04697).

[2]Department of Data Science and Analytics, BI Norwegian Business School, N-0484 Oslo, Norway. Email: fred.e.benth@bi.no. Supported by the SURE-AI Centre grant 357482, Research Council of Norway.

[3]Department of Mathematics, University of Ljubljana, Ljubljana, Slovenia. Email: jan.rems@fmf.uni-lj.si. Supported by the Slovenian Research and Innovation Agency, research core funding No. P1-0448.

# 1 Introduction

SPDEs describe systems in which the spatial dynamics interact with random forcing. They appear in filtering, stochastic reaction diffusion models, turbulent transport, random media, and financial and physical models driven by space–time noise. From a computational point of view, the main difficulty is that the solution

$$U = U(t, x, \omega)$$

depends not only on time and space $(t, x) \in [0, T] \times D$ where $D \subset \mathbb{R}^d$ is a bounded Lipschitz domain, but also on an infinite dimensional variable $\omega$ in the probability space. Directly learning $U$ as a random field is therefore a high dimensional stochastic approximation problem.

The central computational question is therefore: *what should the neural network learn?* Sample wise approaches place the randomness directly in the target and approximate stochastic trajectories or fields. NeSTR instead learns an analytic generating function. More precisely, it approximates the restricted $S$-Transform

$$(t, x, a) \longmapsto u^{(N)}(t, x; a),$$

whose Taylor coefficients at $a = 0$ are exactly the retained Wiener chaos coefficients. Once this deterministic analytic map has been learned, the random field is reconstructed by combining its derivatives with fixed Hermite random variables. Learning and randomness are thereby decoupled: the learning stage produces equation-dependent deterministic coefficient fields, whereas randomness enters through a fixed Hermite basis that can be generated and stored independently. This motivates learning the $S$-Transform instead of stochastic trajectories.

The objective of the NeSTR approach is to replace this stochastic learning problem by a deterministic transform learning problem. Instead of approximating the solution realizations $U(t, x, \omega)$ directly, we approximate the $S$-Transform

$$u(t, x; h) := (SU(t, x))(h),$$

where $h$ is a test function in the Gaussian noise space. For each fixed $h$, the quantity $u(t, x; h)$ is deterministic. Moreover, when the SPDE is written in Wick form, the $S$-Transform converts Wick products into ordinary products. Thus, the transformed object often satisfies a deterministic PDE, or a family of deterministic PDEs parametrized by $h$.

The numerical contribution of NeSTR is to learn this deterministic $S$-Transform map on a finite dimensional restriction of the transform space. After choosing noise modes $e_1, \ldots, e_N$, we write

$$h^{(N)}(a) = \sum_{k=1}^{N} a_k e_k, \qquad a = (a_1, \ldots, a_N),$$

and learn

$$(t, x, a) \longmapsto u^{(N)}(t, x; a) := (SU(t, x))(h^{(N)}(a)).$$

The stochastic solution is then reconstructed by an approximate inverse $S$-Transform. The key identity is that the Wiener chaos coefficients are the Taylor coefficients of the $S$-Transform at the origin:

$$c_\alpha(t, x) = \frac{1}{\alpha!} \partial_a^\alpha u^{(N)}(t, x; 0).$$

Thus the map learned by NeSTR is not merely a surrogate for a deterministic PDE; it is a mechanism for recovering the stochastic expansion of the SPDE solution. This places an unusual demand on the approximation: it is the derivatives of the learned map at a single point, not its values, that carry the answer, and a small training loss does not by itself control them. We therefore do not differentiate a generic network. The dependence on $a$ is represented by a polynomial of fixed degree whose coefficients are neural network fields in $(t, x)$, so that $\partial_a^\alpha u_\theta(t, x; 0)$ is an exact linear readout of learned quantities rather than the output of repeated automatic differentiation.

This separation also gives a re-usability property. The learning stage produces the deterministic fields $c_\alpha(t, x)$, whereas the stochastic stage uses the Hermite random variables

$$H_\alpha(\xi), \qquad \xi = (\xi_k)_{k\geq 1}, \qquad \xi_k := W(e_k),$$

associated with a fixed orthonormal basis $(e_k)_{k\geq 1}$ of the Gaussian noise space. Thus $\xi_k$ are the Gaussian coordinates of the driving noise, and are independent standard normal random variables.

For a fixed Gaussian noise basis and chaos truncation, these Hermite variables can be simulated and stored once. Combining the same stored basis with different equation-dependent coefficient fields then produces samples from different SPDE approximations. Proposition 3.9 formalizes this observation as a transfer principle: for SPDEs driven by the same Gaussian noise, changing the equation changes only the deterministic chaos coefficients, while the retained random Hermite basis remains fixed.

In the literature, one may distinguish between two different paths of the neural PDE literature. The works of Han, Jentzen, E, Beck, and co-authors primarily develop numerical methods for high dimensional deterministic PDEs, often through nonlinear Feynman–Kac or backward stochastic differential equation (BSDE) representations [12, 3, 2]. In NeSTR, these ideas concern the deterministic solver used *after* applying the $S$-Transform. The SPDE literature addresses a different question: how the random field itself is represented, learned, and reconstructed. The relevant comparisons there are with samplewise solvers, modal representations, and Wiener chaos methods.

This viewpoint is related to, but distinct from, other neural approaches for stochastic dynamics and PDEs. Physics-informed neural networks (PINNs) and deep BSDE solvers approximate deterministic PDE solutions or deterministic reformulations of PDEs [25, 12, 3, 2]. Recent latent variable approaches, such as [29], combine spectral Galerkin projection with truncated Wiener chaos and infer latent temporal dynamics from solution observations. In contrast, NeSTR learns the deterministic $S$-Transform itself and uses the inverse $S$-Transform through chaos derivatives for reconstruction. The role of the neural network is therefore to approximate a deterministic analytic map whose local derivatives encode the law of the stochastic solution.

The comparison with existing methods can be summarized as follows. PINNs and deterministic neural PDE solvers learn a PDE solution or residual, but they do not directly reconstruct the Wiener chaos of a stochastic field. The same applies to the two other deterministic families: deep BSDE schemes represent a value function through a backward equation and return a value or gradient approximation, while neural operators learn a solution operator between function spaces and return a deterministic field for each supplied input. Modal space PINN methods combine physics-informed learning with dynamically orthogonal or biorthogonal stochastic representations [30]. Deep learning algorithms for SPDEs have also been tested on stochastic heat equations with additive and multiplicative noise by using a network for each realization of the driving noise [1]. More recently, neural networks have been used directly inside a truncated Wiener chaos expansion of the SPDE solution, with examples including the stochastic heat equation [23]. Latent Wiener chaos methods instead keep a finite chaos basis and learn latent coefficient dynamics from stochastic solution observations [29]. Our proposed NeSTR approach is different: it learns the restricted $S$-Transform $u^{(N)}(t, x; a)$, treats the noise through deterministic transform coordinates $a$, and reconstructs the stochastic field by analytic derivatives at the origin. Thus, the comparison is not only in numerical accuracy, but also in what object is learned: samplewise solution values, modal coefficients, latent chaos dynamics, chaos coefficients themselves, or the transform whose Taylor coefficients are the chaos coefficients. The stochastic output differs accordingly: none at all for the deterministic solvers unless randomness is supplied as an input, one field or trajectory per realization for samplewise solvers, a truncated chaos reconstruction for the modal and chaos methods, and for NeSTR coefficients obtained by differentiation at the origin followed by a Hermite reconstruction.

In short, NeSTR learns an analytic generating function rather than stochastic trajectories. This distinction underlies the remainder of the paper.

There is also a complementary literature on learning maps between infinite dimensional spaces. Benth, Detering, and Galimberti establish universal approximation results for continuous maps and neural networks from Fréchet spaces into Banach spaces and apply related networks in Hilbert space to option pricing on electricity forward curves [4, 5]. These results are close in spirit because the present problem is also an operator approximation problem, although the Hida test and distribution spaces and the derivative based inverse $S$-Transform that we are concerned with are not covered verbatim by that framework. Galimberti, Kratsios, and Livieri develop causal neural operators for infinite dimensional stochastic dynamical systems, while Kratsios, Livieri, and Schmocker propose a neural chaos architecture for adapted square integrable processes [11, 19]. Di Nunno and Díaz Lozano combine Wiener chaos, neural networks, and an operator viewpoint for BSDE solution maps motivated by dynamic risk measures [10]. Finally, Cuchiero, Schmocker, and Teichmann prove global approximation results on weighted infinite dimensional spaces, beyond approximation only on compact sets, and Schmocker and Teichmann extend this viewpoint to differentiable maps and simultaneous approximation of their derivatives [7, 26]. The latter perspective on derivatives is particularly relevant here because the inverse $S$-Transform is recovered from derivatives at the origin.

We develop the NeSTR framework through finite dimensional restriction and restricted $S$-Transforms. We then formulate the deterministic parametric PDEs, learning objectives, inverse reconstruction, and consistency theory. We also explain why continuity of the $S$-Transform and of its inverse on the image justifies learning in transform space. Finally, we give numerical tests that separate the roles of solver accuracy, transform accuracy, and inverse reconstruction. A deterministic benchmark with a closed form solution isolates the underlying solver before any transform is involved. The additive noise heat equation provides a controlled affine case with known chaos coefficients, while the multiplicative Wick heat equation is the natural nonlinear test for recovering higher order chaos. In the multiplicative case the polynomial representation recovers the chaos coefficients more accurately than a multilayer perceptron differentiated automatically, with the gap widening as the order increases. We also show that the reconstruction error splits exactly into a recovery part and a chaos truncation part, which identifies the usable truncation order and, at the same time, the current limitation of the method: beyond that order the exact coefficients fall below the attainable accuracy.

The remainder of the paper is organized as follows. Section 2 recalls the white noise tools used in the sequel. Section 3 is the main methodological section and develops the NeSTR framework. Section 4 presents the numerical approximation scheme, implementation details, and neural experiments. Section 5 concludes and summarizes the distinctive contribution and limitations of NeSTR.

## 2 Preliminaries on White Noise Tools

We briefly recall the standard white noise notation needed for the NeSTR framework. The construction of the Hida test space, the Hida distribution space, Wiener chaos, Wick products, and the $S$-Transform are classical; see, for example, [13, 14, 20, 15]. We use only the facts listed below.

Let $\mathcal{H}$ be a real separable Hilbert space, typically $\mathcal{H} = L^2([0,T] \times D)$, and let $W = \{W(h) : h \in \mathcal{H}\}$ be an isonormal Gaussian process: Two levels of structure should be distinguished here. The square integrable Wiener chaos construction requires only the Gaussian Hilbert space $\mathcal{H}$. To introduce Hida test functions and distributions, one additionally fixes a nuclear Gelfand triple

$$E \subset \mathcal{H} \subset E^*,$$

with dense continuous embeddings, and constructs the corresponding Hida triple over $E$. For a temporally driven SPDE one may take $\mathcal{H} = L^2([0,T])$, in which case the spatial point $x \in D$ is an external parameter of the chaos coefficients. For space–time white noise, as in the equations containing $\dot{W}(t,x)$ below, the natural Gaussian space is $\mathcal{H}$, so both temporal and spatial noise directions enter the chaos expansion. Standard constructions for these SPDE settings are given in [9, 15]. One has,

$$\mathbb{E}[W(h)W(g)] = \langle h, g \rangle_{\mathcal{H}} .$$

Fix an orthonormal basis $(e_k)_{k\geq 1}$ of $\mathcal{H}$ and set

$$\xi_k := W(e_k), \qquad k \geq 1.$$

Then $(\xi_k)_{k\geq 1}$ are independent standard normal variables and, for $h = \sum_k h_k e_k$,

$$W(h) = \sum_{k\geq 1} h_k \xi_k \quad \text{in } L^2(\Omega).$$

Every $F \in L^2(\Omega)$ admits a Wiener Itô chaos expansion [6, 17]

$$F = \sum_{n=0}^{\infty} I_n(f_n), \qquad f_n \in \mathcal{H}^{\odot n},$$

Here $\mathcal{H}^{\odot n}$ denotes the $n$ fold symmetric Hilbert tensor product of $\mathcal{H}$, with the convention $\mathcal{H}^{\odot 0} = \mathbb{R}$. Furthermore, it holds that,

$$\|F\|^2_{L^2(\Omega)} = \sum_{n=0}^{\infty} n! \; \|f_n\|^2_{\mathcal{H}^{\odot n}} .$$

Let $\mathcal{J}$ denote the set of finite support multi-indices,

$$\mathcal{J} := \{\alpha = (\alpha_k)_{k\geq 1} : \alpha_k \in \mathbb{N}_0, \; \#\operatorname{supp}(\alpha) < \infty\},$$

and write

$$|\alpha| := \sum_{k\geq 1} \alpha_k, \qquad \alpha! := \prod_{k\geq 1} \alpha_k!.$$

Equivalently, using the probabilists' Hermite polynomials, the chaos expansion can be written as

$$F = \sum_{\alpha \in \mathcal{J}} c_\alpha H_\alpha(\xi), \qquad H_\alpha(\xi) := \prod_{k\geq 1} H_{\alpha_k}(\xi_k),$$

where

$$\|F\|_{L^2(\Omega)}^2 = \sum_{\alpha \in \mathcal{J}} \alpha! \, |c_\alpha|^2.$$

To include generalized stochastic objects we use the Hida triple

$$(\mathcal{S}) \subset L^2(\Omega) \subset (\mathcal{S})^*.$$

Thus a solution may be interpreted either as a square integrable random variable or, when necessary, as a Hida distribution.

The Wick product $\diamond$ is the product compatible with chaos. In Hermite coordinates it satisfies

$$H_\alpha(\xi) \diamond H_\beta(\xi) = H_{\alpha+\beta}(\xi).$$

The Wick exponential of a linear Gaussian variable is

$$\exp^\diamond(W(h)) = \exp\left(W(h) - \frac{1}{2}\|h\|_H^2\right).$$

For $F \in L^2(\Omega)$, the $S$-Transform is

$$(SF)(h) := \mathbb{E}\left[F \exp^\diamond(W(h))\right], \qquad h \in \mathcal{H}_\mathbb{C},$$

and it extends to $(\mathcal{S})^*$ by duality on the complexified Hida test space $E_\mathbb{C}$. Thus the natural transform domain is $\mathcal{H}_\mathbb{C}$ for square integrable random variables and $E_\mathbb{C}$ for general Hida distributions. The following standard properties are the only ones used in the sequel; see [14, 20, 15] for the full Hida space statement.

The next theorem is classical. We include its short proof because the paper is formulated for an arbitrary real separable Gaussian Hilbert space $\mathcal{H}$, rather than only the canonical realization $\mathcal{H} = L^2(\mathbb{R})$. At the square integrable level the proof is basis independent and transfers directly to any such $\mathcal{H}$. No new $S$-Transform theorem is claimed; the genuinely distributional extension uses the nuclear triple above and is taken from the standard white noise literature.

**Theorem 2.1** (Standard $S$-Transform facts)**.** *Let $E \subset \mathcal{H} \subset E^*$ be the nuclear triple underlying the Hida spaces. For $F, G \in (\mathcal{S})^*$, the following hold on $E_\mathbb{C}$; in the square integrable case the $S$-Transform extends to $\mathcal{H}_\mathbb{C}$:*

*(i) $h \mapsto (SF)(h)$ is entire on $E_\mathbb{C}$, and on $\mathcal{H}_\mathbb{C}$ when $F \in L^2(\Omega)$;*

*(ii) the $S$-Transform is injective on $(\mathcal{S})^*$;*

*(iii) whenever the Wick product is defined,*

$$S(F \diamond G) = (SF)(SG); \tag{1}$$

*(iv) for $g \in \mathcal{H}$,*

$$S(W(g))(h) = \langle g, h \rangle_\mathcal{H}. \tag{2}$$

*Proof.* We prove the square integrable statement on an arbitrary separable $\mathcal{H}$, which contains the identities used later in the paper. The extension to Hida distributions is not reproved: it follows from the standard duality construction and the characterization theorem for $S$-Transforms on the nuclear test space; see the references above.

Let $F \in L^2(\Omega)$ have Hermite chaos expansion

$$F = \sum_\alpha c_\alpha H_\alpha(\xi), \qquad \sum_\alpha \alpha! \, |c_\alpha|^2 < \infty.$$

For $h = \sum_k z_k e_k$, the Wick exponential factorizes as

$$\exp^\diamond(W(h)) = \prod_{k \geq 1} \exp^\diamond(z_k \xi_k) = \prod_{k \geq 1} \sum_{m=0}^{\infty} \frac{z_k^m}{m!} H_m(\xi_k).$$

Thus, in multi-index notation,

$$\exp^\diamond(W(h)) = \sum_\beta \frac{z^\beta}{\beta!} H_\beta(\xi),$$

with convergence in the usual exponential vector sense. Using $\mathbb{E}[H_\alpha H_\beta] = \alpha! \delta_{\alpha\beta}$, we obtain

$$(SF)(h) = \mathbb{E}\left[\sum_\alpha c_\alpha H_\alpha(\xi) \sum_\beta \frac{z^\beta}{\beta!} H_\beta(\xi)\right] = \sum_\alpha c_\alpha z^\alpha.$$

This is the coordinate $\mathcal{H}$-Transform of $F$. On every finite dimensional restriction it is a power series in finitely many complex variables, and the standard $S$-Transform estimates imply the corresponding entire extension on $\mathcal{H}_\mathbb{C}$. This proves analyticity.

The same formula proves injectivity in $L^2$. If $(SF)(h) = 0$ for all $h$, then all Taylor coefficients at the origin vanish:

$$c_\alpha = \frac{1}{\alpha!} \partial_z^\alpha (SF)(0) = 0.$$

Hence every chaos coefficient of $F$ is zero, and therefore $F = 0$ in $L^2(\Omega)$. The Hida distribution injectivity is obtained by the same coefficient argument in the Hida chaos coordinates.

For the Wick product identity, it is enough to check chaos basis elements and then extend by bilinearity and continuity whenever the Wick product is defined. Since

$$H_\alpha(\xi) \diamond H_\beta(\xi) = H_{\alpha+\beta}(\xi),$$

we have

$$S(H_\alpha \diamond H_\beta)(h) = S(H_{\alpha+\beta})(h) = z^{\alpha+\beta} = z^\alpha z^\beta = (SH_\alpha)(h)(SH_\beta)(h).$$

Therefore $S(F \diamond G) = (SF)(SG)$ for finite chaos polynomials and, by the standard definition of the Wick product, for all pairs for which the product is defined.

Finally, take $F = W(g)$, where $g = \sum_k g_k e_k$. Then $W(g) = \sum_k g_k \xi_k$, so only first order chaos terms are present. The previous coefficient formula gives

$$S(W(g))(h) = \sum_k g_k z_k = \langle g, h \rangle_H \,,$$

which is (2). This completes the proof of the identities used in the sequel. □

**Theorem 2.2** (Potthoff–Streit characterization and inverse continuity)**.** *Let $\mathcal{U}$ denote the Potthoff–Streit space of U-functionals, namely the entire functions $\Phi$ on the complexified white noise test space such that, for suitable constants $C, K > 0$ and a Hilbertian test function norm $\| \cdot \|_p$,*

$$|\Phi(z\varphi)| \leq C \exp\bigl(K|z|^2 \|\varphi\|_p^2\bigr), \qquad z \in \mathbb{C},$$

*with the corresponding standard locally convex topology of uniform control by such exponential bounds. Then the $S$-Transform is a topological isomorphism*

$$S : (\mathcal{S})^* \longrightarrow \mathcal{U}.$$

*In particular, $S$ is continuous and injective, and its inverse*

$$S^{-1} : \mathcal{U} \longrightarrow (\mathcal{S})^*$$

*is continuous.*

This is the quoted Potthoff–Streit characterization theorem for $S$-Transforms of Hida distributions [24]. We do not prove this Hida space result here; it is used only as a functional analytic stability statement for the inverse transform. All explicit coefficient calculations in this paper are carried out in the square integrable chaos setting. SPDE oriented accounts of the white noise setting can be found in [9, 15].

*Remark* 2.3. It is well-known (see e.g. Walsh [28]) that stochastic partial differential equations may be singular in the sense of not possessing any solution being a random variable. Holden *et al.* [15] show that many equations of interest possess solutions which are Hida distributions, or even more generally, Kondratiev distributions. One may additionally construct appropriate spaces of Hida (or Kondratiev distributions) which also control the regularity in time and space via Sobolev or generalised function spaces, say, for the solutions. In this paper we are focused on parabolic stochastic partial differential equations which are analysed in detail in Holden *et al.* [15], and known to have unique solutions in appropriate spaces. As we are concerned with developing a deep learning approach to solving such equations, we will not discuss existence and uniqueness in more detail here, nor present the broad theory for spaces where the solutions live, but refer to reader to Holden *et al.* [15].

As already indicated in the proof of Theorem 2.1 above, the coordinate form of the $S$-Transform is the $\mathcal{H}$-Transform. If $h = \sum_k z_k e_k$, then

$$(\mathcal{H}F)(z) = (SF)(h).$$

Consequently, if $(\mathcal{H}F)(z) = \sum_\alpha c_\alpha z^\alpha$, then

$$c_\alpha = \frac{1}{\alpha!} \partial_z^\alpha (\mathcal{H}F)(0).$$

This derivative formula is the inverse reconstruction identity exploited in Section 3.

### 2.1 Example: $S$-Transform of SPDEs

This subsection illustrates the abstract transform identities on SPDEs. The goal is to show, before introducing the learning procedure, how Wick stochastic terms become deterministic coefficients after applying the $S$-Transform.

Let $U = U(t, x)$ be a stochastic field. We define its transform by

$$u(t, x; h) := (SU(t, x))(h), \qquad h \in E_{\mathbb{C}},$$

with the larger domain $h \in \mathcal{H}_{\mathbb{C}}$ whenever $U(t, x) \in L^2(\Omega)$. By injectivity, $u(t, x; \cdot)$ determines $U(t, x)$. If the SPDE is written in Wick form, then (1) turns Wick products into ordinary products and replaces stochastic forcing by deterministic dependence on $h$. This is why the $S$-Transform often converts an SPDE into a deterministic parametric PDE.

For example, consider the additive noise heat equation on $D$,

$$\begin{cases} \partial_t U(t,x) = \Delta U(t,x) + \dot{W}(t,x), \\ U(0,x) = u_0(x), \\ U(t,x) = 0 \quad \text{on } \partial D, \end{cases} \tag{3}$$

where the additive stochastic forcing is interpreted in the Wick, equivalently Skorohod, sense and the initial condition $u_0$ is a sufficiently regular deterministic function.

Let $(G_t)_{t\geq 0}$ denote the Dirichlet heat kernel on $D$. The corresponding mild solution is

$$U(t,x) = \int_D G_t(x,y)u_0(y)\,dy + \int_0^t \int_D G_{t-s}(x,y)\,W(ds,dy).$$

Since the integrand in the stochastic convolution is deterministic, the Skorohod integral coincides in this case with the usual Itô–Walsh stochastic integral. We refer to [15] for more details. Applying the $S$-Transform gives

$$\begin{cases} \partial_t u(t,x;h) = \Delta u(t,x;h) + h(t,x), \\ u(0,x;h) = u_0(x), \\ u(t,x;h) = 0 \quad \text{on } \partial D. \end{cases} \tag{4}$$

Equivalently, with the Dirichlet heat kernel $G_t$,

$$u(t,x;h) = \int_D G_t(x,y)u_0(y)\,dy + \int_0^t \int_D G_{t-s}(x,y)h(s,y)\,dy\,ds. \tag{5}$$

For fixed $h$, this deterministic PDE also admits the Feynman–Kac form [18]

$$u(t,x;h) = \mathbb{E}\left[u_0(X_t^x)\mathbb{1}_{\{t<\tau_D\}} + \int_0^{t\wedge\tau_D} h(t-s, X_s^x)\,ds\right]. \tag{6}$$

Here $X^x$ is the diffusion (in fact, a Brownian motion in $\mathbb{R}^d$ with speed $\sqrt{2}$) started at $x$ with infinitesimal generator $\Delta$, and $\tau_D$ is its first exit time from $D$. Thus the learning target is deterministic in $(t,x,h)$, not random in $\omega$.

## 3 The NeSTR Framework

This is the main section of the paper. Its purpose is to define NeSTR by turning the abstract white noise transform into a concrete numerical procedure. The starting point is a stochastic field $U(t,x)$, but the learning target will not be $U$ itself. Instead, we learn the deterministic $S$-Transform map

$$h \longmapsto (SU(t,x))(h).$$

This is useful for two reasons. First, the $S$-Transform removes the explicit random component $\omega$: after transformation, the unknown is a deterministic function of $(t,x,h)$. Second, the Taylor coefficients of this deterministic function at the origin are exactly the Wiener chaos coefficients of $U(t,x)$. Thus the NeSTR approximation of the transform map provides a direct route to reconstructing the stochastic solution. Or, put equivalently, NeSTR learns an analytic generating function rather than stochastic trajectories; its derivatives encode the stochastic solution.

At the square integrable level the full variable belongs to $\mathcal{H}_{\mathbb{C}}$; for a general Hida distribution it belongs to $E_{\mathbb{C}}$. In either case the transform variable is infinite dimensional, so a computable method must restrict it to a finite dimensional subspace. We choose a basis $(e_k)_{k\geq 1}$ of the noise space and write

$$h^{(N)}(a) = \sum_{k=1}^{N} a_k e_k, \qquad a = (a_1, \dots, a_N) \in \mathbb{C}^N.$$

The object learned by NeSTR is therefore

$$u^{(N)}(t, x; a) := (SU(t, x))(h^{(N)}(a)).$$

Hence, NeSTR learns the restricted solution operator

$$\mathcal{G}_N : a \longmapsto u^{(N)}(\cdot, \cdot; a).$$

This may be viewed as an operator learning problem as the output is an entire space–time field and the input is a function, although represented in its coordinate weights along basis functions. DeepONets and Fourier neural operators provide architectures for operator learning problems [22, 21]. Universal approximation results in general infinite dimensional spaces and causal operator settings are developed in [4, 11]; NeSTR adds the specific analytic structure that derivatives of $\mathcal{G}_N$ at the origin recover Wiener chaos coefficients.

The stochastic approximation is then obtained by differentiating the learned map with respect to $a$ at the origin and inserting the resulting coefficients into a truncated Hermite expansion.

In other words, NeSTR has four steps:

$$U(t, x) \xrightarrow{S} u^{(N)}(t, x; a) \xrightarrow{\text{learn}} u_\theta(t, x; a) \xrightarrow{\partial_a^\alpha|_{a=0}} \widehat{c}_\alpha(t, x) \xrightarrow{\text{Hermite series}} \widehat{U}_{N,M}(t, x).$$

The key reconstruction identity is

$$c_\alpha(t, x) = \frac{1}{\alpha!} \partial_a^\alpha u^{(N)}(t, x; 0).$$

The subsections below describe each part of the NeSTR pipeline: finite dimensional restriction, restricted $S$-Transform, transformed deterministic PDE, learning objective, training data, inverse reconstruction, and consistency.

### 3.1 Finite dimensional restriction and restricted transform

This subsection defines the finite dimensional transform coordinates used by NeSTR. It explains how the infinite dimensional test function $h$ is replaced by finitely many coefficients $a_1, \dots, a_N$.

Recall that $\mathcal{H}$ is a real separable Hilbert space, with $(e_k)_{k\geq 1}$ a fixed orthonormal basis of $\mathcal{H}$. When the solution is treated only as a Hida distribution, we choose the basis with $e_k \in E$; such a basis can be constructed by applying Gram–Schmidt to a countable dense subset of $E$. In the SPDE setting of interest, the canonical example is

$$\mathcal{H} = L^2([0, T] \times D),$$

where $D \subset \mathbb{R}^d$ is a bounded spatial domain.

For $N \in \mathbb{N}$, define the finite dimensional subspace

$$\mathcal{H}^{(N)} := \operatorname{span}\{e_1, \dots, e_N\} \subset \mathcal{H}.$$

Its complexification is

$$\mathcal{H}_{\mathbb{C}}^{(N)} := \mathrm{span}_{\mathbb{C}}\{e_1, \dots, e_N\} \subset \mathcal{H}_{\mathbb{C}}.$$

For every $a = (a_1, \dots, a_N) \in \mathbb{C}^N$, define

$$h^{(N)}(a) := \sum_{k=1}^{N} a_k e_k \in \mathcal{H}_{\mathbb{C}}^{(N)}. \tag{7}$$

Thus the infinite dimensional parameter $h \in \mathcal{H}_{\mathbb{C}}$ is replaced by the finite dimensional coordinate vector $a \in \mathbb{C}^N$. The integer $N$ determines how many directions of the driving noise are retained. Increasing $N$ gives a richer approximation of the random input, while keeping $N$ finite makes the learning problem finite dimensional.

This restriction should be interpreted as a projection in transform space, not as a discretization of the variables $(t, x)$. The transformed PDE may still be solved or sampled in continuous space, while the random direction is represented by the coordinates $a_1, \dots, a_N$. This separation is important: NeSTR controls stochastic complexity through $N$, and physical space complexity through the PDE solver or neural architecture.

The choice of basis determines which features of the noise are retained first. For space–time white noise, one may use Fourier modes, eigenfunctions of a spatial operator, wavelets, or any orthonormal system adapted to the domain. The NeSTR itself does not depend on this choice. What matters is that the first $N$ modes provide a finite coordinate system in which the transform can be sampled, learned, and differentiated.

*Remark* 3.1. The map $a \mapsto h^{(N)}(a)$ is a linear isomorphism from $\mathbb{C}^N$ onto $\mathcal{H}_{\mathbb{C}}^{(N)}$. Therefore one may equivalently formulate the finite dimensional transform problem in terms of $a \in \mathbb{C}^N$ or $h^{(N)} \in \mathcal{H}_{\mathbb{C}}^{(N)}$.

**Restricted $S$-Transform of a stochastic field.** Let

$$U = U(t, x)$$

be a stochastic field defined on $[0, T] \times D$ with values in $(\mathcal{S})^*$, or more specifically in $L^2(\Omega)$ whenever pointwise square integrability is available. For each fixed $(t, x)$, the map

$$h \longmapsto u(t, x; h) = (SU(t, x))(h)$$

is, by Theorem 2.1, defined on $E_{\mathbb{C}}$ in the Hida distribution setting. If $U(t, x) \in L^2(\Omega)$, the transform extends to $\mathcal{H}_{\mathbb{C}}$.

**Definition 3.2** (Restricted $S$-Transform map)**.** For $N \in \mathbb{N}$, the *restricted $S$-Transform* of $U$ is the map

$$u^{(N)} : [0, T] \times D \times \mathbb{C}^N \longrightarrow \mathbb{C}, \qquad u^{(N)}(t, x; a) := u(t, x; h^{(N)}(a)). \tag{8}$$

For a general Hida distribution this definition presupposes $e_1, \dots, e_N \in E$, so that $h^{(N)}(a) \in E_{\mathbb{C}}$.

For each fixed $(t, x)$, the function $a \mapsto u^{(N)}(t, x; a)$ is simply the restriction of the entire function $h \mapsto u(t, x; h)$ to the finite dimensional subspace $\mathcal{H}_{\mathbb{C}}^{(N)}$. Thus, as we have stated previously, the neural network therefore does not approximate a random sample path. It approximates a deterministic analytic function whose Taylor expansion encodes the stochastic structure of $U(t, x)$.

This distinction is central. If one learns sample paths directly, the network must represent both deterministic evolution and random fluctuations at the same time. In the transform formulation, the random fluctuations are encoded by the coordinate variable $a$, and the object to be learned is a deterministic response surface. The inverse step then translates this response surface back into Wiener chaos coefficients.

**Proposition 3.3** (Analyticity of the restricted transform)**.** *Let $(t,x) \in [0,T] \times D$ be fixed. Assume either that $U(t,x) \in L^2(\Omega)$, or that $U(t,x) \in (\mathcal{S})^*$ and the retained modes $e_1, \dots, e_N$ belong to $E$. Then the map*

$$a \longmapsto u^{(N)}(t,x;a)$$

*is entire on $\mathbb{C}^N$.*

*Proof.* Fix $(t,x)$ and write

$$F := U(t,x) \in (\mathcal{S})^*.$$

By Theorem 2.1, the $S$-Transform $\mathcal{S}_F(h) := (SF)(h)$ is entire on $\mathcal{H}_\mathbb{C}$ in the square integrable case and on $E_\mathbb{C}$ in the Hida distribution case. Let $\mathcal{X}_\mathbb{C}$ denote the corresponding transform domain.

Define the finite dimensional restriction map

$$L_N : \mathbb{C}^N \longrightarrow \mathcal{X}_\mathbb{C}, \qquad L_N(a) := h^{(N)}(a) = \sum_{j=1}^N a_j e_j.$$

In the Hida distribution case its range is contained in $E_\mathbb{C}$ by the assumption on the retained modes.

Since $\mathcal{S}_F$ is entire at the origin, it has a convergent power series in the noise directions. Restricting this series to the span of $e_1, \dots, e_N$ gives

$$u^{(N)}(t,x;a) = (SF)\left(\sum_{j=1}^N a_j e_j\right) = \sum_{\alpha \in \mathbb{N}_0^N} c_\alpha(t,x) a^\alpha$$

locally in $a$, where the coefficients are precisely the derivatives of the $S$-Transform in the retained directions. Because the original $S$-Transform is entire, this restricted power series has infinite radius of convergence in every finite dimensional direction. Hence the restricted map is an entire function on $\mathbb{C}^N$. □

*Remark* 3.4. Proposition 3.3 is one of the key structural facts in NeSTR. The object learned numerically is not just an arbitrary parametric map, but a finite dimensional restriction of an entire function.

### 3.2 Additive noise heat equation

We return to the additive noise heat equation of Example 2.1, whose $S$-Transform is given by (5). Here our purpose is different: we now restrict the transform variable to the finite dimensional family $h^{(N)}(a)$ defined in (7) and identify the resulting dependence on the transform coordinates $a$. Substituting $h = h^{(N)}(a)$ from (7) into (5) yields the restricted transform

$$u^{(N)}(t,x;a) = \int_D G_t(x,y) u_0(y)\, dy + \sum_{k=1}^N a_k \int_0^t \int_D G_{t-s}(x,y) e_k(s,y)\, dy\, ds. \tag{9}$$

Therefore, for the additive noise heat equation, the restricted transform depends *linearly* on the finite dimensional parameter $a$.

**Proposition 3.5** (Continuity and linear growth in the additive case)**.** *For each fixed $(t,x) \in [0,T] \times D$, the map*

$$a \longmapsto u^{(N)}(t,x;a)$$

*is affine linear and therefore Lipschitz continuous on $\mathbb{C}^N$. More precisely,*

$$u^{(N)}(t,x;a) = c_0(t,x) + \sum_{k=1}^{N} a_k c_k(t,x),$$

*where*

$$c_0(t,x) := \int_D G_t(x,y)u_0(y)\,dy, \qquad c_k(t,x) := \int_0^t \int_D G_{t-s}(x,y)e_k(s,y)\,dy\,ds.$$

*Proof.* Fix $(t,x) \in [0,T] \times D$. Starting from the $S$-Transform representation (9), we now restrict the transform variable to

$$h = h^{(N)}(a) = \sum_{k=1}^{N} a_k e_k, \quad e_1,\dots,e_N, \qquad a = (a_1,\dots,a_N) \in \mathbb{C}^N.$$

Substituting this expression into the preceding formula gives

$$\begin{aligned} u^{(N)}(t,x;a) &= u(t,x;h^{(N)}(a)) \\ &= \int_D G_t(x,y)u_0(y)\,dy + \int_0^t \int_D G_{t-s}(x,y)\left(\sum_{k=1}^{N} a_k e_k(s,y)\right)dy\,ds. \end{aligned}$$

Because the sum over $k$ is finite, it can be moved outside the integral without any limiting argument. Hence

$$u^{(N)}(t,x;a) = \int_D G_t(x,y)u_0(y)\,dy + \sum_{k=1}^{N} a_k \int_0^t \int_D G_{t-s}(x,y)e_k(s,y)\,dy\,ds.$$

Define

$$c_0(t,x) := \int_D G_t(x,y)u_0(y)\,dy$$

and, for $k = 1,\dots,N$,

$$c_k(t,x) := \int_0^t \int_D G_{t-s}(x,y)e_k(s,y)\,dy\,ds.$$

Then

$$u^{(N)}(t,x;a) = c_0(t,x) + \sum_{k=1}^{N} a_k c_k(t,x),$$

which proves the affine representation.

It remains only to spell out continuity and linear growth. Since $(t,x)$ is fixed, the numbers $c_0(t,x), c_1(t,x),\dots,c_N(t,x)$ are constants with respect to $a$. Therefore, for $a, b \in \mathbb{C}^N$,

$$|u^{(N)}(t,x;a) - u^{(N)}(t,x;b)| = \left|\sum_{k=1}^{N}(a_k - b_k)c_k(t,x)\right|.$$

By the Cauchy–Schwarz inequality,

$$|u^{(N)}(t,x;a) - u^{(N)}(t,x;b)| \le \left(\sum_{k=1}^{N}|c_k(t,x)|^2\right)^{1/2} |a-b|_{\mathbb{C}^N}.$$

Thus the map is Lipschitz, and in particular continuous. Taking $b = 0$ gives

$$|u^{(N)}(t,x;a)| \le |c_0(t,x)| + \left(\sum_{k=1}^{N}|c_k(t,x)|^2\right)^{1/2} |a|_{\mathbb{C}^N}.$$

This is the claimed linear growth in the transform variable $a$. The proof is complete. □

In this example the inverse reconstruction is especially transparent. Since $u^{(N)}$ is affine in $a$, all second and higher derivatives with respect to $a$ vanish. Therefore the only nonzero chaos coefficients in the retained modes are the deterministic coefficient $c_0$ and the first order coefficients $c_1, \dots, c_N$. This makes the additive equation a useful diagnostic test: any learned curvature in $a$ corresponds to spurious higher order chaos.

### 3.3 Deterministic learning problem and training data

This subsection turns the transform formulation into a supervised deterministic learning problem. It identifies the neural target, the training region, and the role of samples near the origin in transform space.

The additive heat equation shows the mechanism in a case where all formulas are explicit. We now describe the same reduction at the level needed for learning. In more general SPDEs, after applying the $S$-Transform one typically obtains a deterministic parametric PDE of the form

$$\mathcal{A}(t, x, u(\cdot, \cdot; h), h) = 0,$$

with suitable initial and boundary conditions, where $\mathcal{A}$ is a deterministic operator depending on the original SPDE. Restricting to $\mathcal{H}_{\mathbb{C}}^{(N)}$ yields the finite dimensional parametric PDE

$$\mathcal{A}(t, x, u^{(N)}(\cdot, \cdot; a), h^{(N)}(a)) = 0, \qquad a \in \mathbb{C}^N. \tag{10}$$

The practical meaning of this reduction is that the stochastic problem has been transformed into a family of deterministic PDEs indexed by a finite dimensional parameter $a \in \mathbb{C}^N$. For each value of $a$, one solves or approximates an ordinary deterministic problem. The parameter $a$ records how the solution responds to perturbations in the retained noise directions. This response curve is precisely what must be learned, because its derivatives at $a = 0$ produce the chaos coefficients.

Thus NeSTR changes the role of randomness. Randomness is no longer represented by repeated sampling of $\omega$. Instead, it appears as a set of deterministic coordinates in the transform variable. This is advantageous for machine learning: the network sees a supervised deterministic regression problem, while the stochastic interpretation is recovered analytically through the inverse transform.

**Learning objective.** The NeSTR learning problem is now clear. For fixed $N$, one wants to approximate the deterministic map

$$(t, x, a) \longmapsto u^{(N)}(t, x; a),$$

where $a \in \mathbb{C}^N$ parametrizes the restriction of the $S$-Transform to the finite dimensional subspace $\mathcal{H}_{\mathbb{C}}^{(N)}$. The input variables have two separate meanings: $(t, x)$ are the physical variables of the SPDE, while $a$ are transform coordinates describing the retained noise directions. The output is the $S$-Transform value, not a simulated realization of the stochastic solution.

On a compact training region

$$K \subset [0, T] \times D \times \mathbb{C}^N,$$

one may train a neural network by minimizing a supervised loss of the form

$$\int_K |u_\theta(t, x; a) - u^{(N)}(t, x; a)|^2 \, d\mu(t, x, a),$$

where $\mu$ is the sampling distribution used to generate training points. In practice the exact integral is replaced by an empirical average over sampled triples $(t, x, a)$. The choice of sampling measure matters because the recovery of chaos coefficients depends on local behavior near $a = 0$, while global accuracy over the training box controls the quality of the transform surrogate.

The use of a compact set is consistent with the usual scope of classical universal approximation theorems: continuous finite dimensional maps can be approximated uniformly on compact sets [8, 16], and operator versions such as DeepONet similarly formulate approximation on compact subsets of the input function space [22]. The universal approximation theorem on Fréchet spaces of Benth, Detering, and Galimberti provides a related infinite dimensional formulation [4], while weighted space results by Cuchiero, Schmocker and Teichmann [7] allow for global approximation beyond compact sets. Here we make only a local training claim on $K$. Moreover, approximation of transform values alone does not automatically control the derivatives at $a = 0$ needed for inverse reconstruction; this is why the later consistency results impose derivative control and the numerical scheme uses polynomial structure and enhanced sampling near the origin.

**Definition 3.6** (NeSTR approximation of the restricted $S$-Transform)**.** Let $N \in \mathbb{N}$. A neural approximation of the restricted $S$-Transform is a function

$$u_\theta : [0, T] \times D \times \mathbb{C}^N \longrightarrow \mathbb{C},$$

depending on trainable parameters $\theta$, such that

$$u_\theta(t, x; a) \approx u^{(N)}(t, x; a) \qquad \text{for } (t, x, a) \text{ in the region of interest.} \tag{11}$$

In practice, one identifies $\mathbb{C}^N$ with $\mathbb{R}^{2N}$ by writing

$$a_k = \alpha_k + i\beta_k, \qquad (\alpha_1, \beta_1, \dots, \alpha_N, \beta_N) \in \mathbb{R}^{2N},$$

and the neural network is implemented as a real-valued network producing the real and imaginary parts of $u^{(N)}$. In the experiments below the relevant solutions are real on real $a$, so one can train on real boxes $a \in [-A, A]^N$. The complex analyticity remains important theoretically, because it justifies recovering chaos coefficients by derivatives at the origin.

*Remark* 3.7. The advantage of NeSTR is that the learning target in Definition 3.6 is not a random field, but a deterministic parametric map. Once the stochastic problem has been transferred into transform space, one may use standard tools for deterministic PDE approximation, including PINNs, Feynman–Kac synthetic data, or operator learning architectures, including DeepONets and Fourier neural operators [25, 12, 22, 21].

**Training data.** The training data are generated in transform space. Thus each training point is of the form

$$(t, x, a, u^{(N)}(t, x; a)),$$

where $(t, x)$ are physical variables and $a$ is the finite dimensional transform coordinate. The important point is that the target value is deterministic. No sample of the original random field $U(t, x, \omega)$ is used as a target.

When the transformed equation has an explicit semigroup representation, the values of $u^{(N)}$ can be evaluated directly. This is the case for the additive noise heat equation, where (9) gives the exact restricted transform. In more general situations one may compute the same quantity by solving the deterministic parametric PDE (10) for different values of $a$. If a Feynman–Kac formula is available, the training values may also be produced by Monte Carlo simulation of the associated deterministic diffusion problem. In all cases, the data are evaluations of the deterministic $S$-Transform map, not realizations of the stochastic solution.

For reconstruction, it is useful to sample the transform variable $a$ so that the network learns accurately near the origin. The values near $a = 0$ determine the derivatives used in the inverse $S$-Transform, whereas values farther from the origin help stabilize the learned analytic dependence. This is one reason for using polynomial or spectral representations in $a$: they make the derivative extraction step more reliable than differentiating an arbitrary black box network many times.

### 3.4 Inverse reconstruction and consistency

This subsection explains how the learned deterministic transform is converted back into a stochastic field. It also records the consistency and convergence statements that justify the truncated inverse reconstruction.

**Recovery of chaos coefficients.** This is the central mathematical step of NeSTR. The restricted $S$-Transform is useful not merely because it is deterministic, but because it contains the coefficients of the stochastic solution. These coefficients are obtained by differentiating the restricted transform at the origin.

Let $F \in L^2(\Omega)$ and define, for $a \in \mathbb{C}^N$,

$$f^{(N)}(a) := (SF)(h^{(N)}(a)).$$

Write the full Hermite expansion as

$$F = \sum_{\beta \in \mathcal{J}} c_\beta H_\beta(\xi).$$

For $\alpha \in \mathbb{N}_0^N$, we use the same symbol for its zero extension to an element of $\mathcal{J}$, and $c_\alpha$ denotes the corresponding coefficient. We have:

**Theorem 3.8** (Inverse $S$-Transform through chaos derivatives)**.** *For every multi-index $\alpha = (\alpha_1, \dots, \alpha_N) \in \mathbb{N}_0^N$, it holds that*

$$c_\alpha = \frac{1}{\alpha!} \partial_a^\alpha f^{(N)}(0).$$

*Proof.* On the finite dimensional subspace generated by $e_1, \dots, e_N$, the $S$-Transform has the Taylor expansion

$$f^{(N)}(a) = (SF)(h^{(N)}(a)) = \sum_{\alpha \in \mathbb{N}_0^N} c_\alpha a^\alpha, \qquad a^\alpha := a_1^{\alpha_1} \cdots a_N^{\alpha_N},$$

where only coefficients supported in $\{1, \dots, N\}$ are detected. Indeed, if a chaos coefficient involves some mode $k > N$, then the corresponding coordinate of $h^{(N)}(a)$ is $z_k = 0$, so its contribution to $(SF)(h^{(N)}(a))$ is zero. Differentiating this power series at $a = 0$ gives

$$\partial_a^\alpha f^{(N)}(0) = \alpha! \, c_\alpha.$$

This proves the formula. □

Applying this pointwise to $F = U(t, x)$ gives

$$c_\alpha(t, x) = \frac{1}{\alpha!} \partial_a^\alpha u^{(N)}(t, x; 0). \tag{12}$$

This identity is the heart of NeSTR. It says that the inverse $S$-Transform is implemented, after finite dimensional restriction, by taking derivatives at the origin. Therefore the network must learn the transform accurately not only as a function value, but also locally near $a = 0$, where the derivatives determine the stochastic expansion. This is why polynomial or analytic representations in $a$ are natural for the recovery step.

**Approximate reconstruction.** Suppose that a neural network $u_\theta$ approximates $u^{(N)}$ as in (11). The learned transform is now used as a surrogate for the exact $S$-Transform. We define the approximate chaos coefficients by

$$\widehat{c}_\alpha(t,x) := \frac{1}{\alpha!}\partial_a^\alpha u_\theta(t,x;0), \qquad \alpha \in \mathbb{N}_0^N. \tag{13}$$

In principle, these derivatives can be computed by automatic differentiation. In practice, high order derivatives of a generic neural network may be unstable. This motivates the polynomial in $a$ representation used later: if the $a$-dependence of $u_\theta$ is represented by an explicit polynomial basis, then the derivatives at $a = 0$ are obtained exactly from the learned polynomial coefficients.

For fixed $N \in \mathbb{N}$ and $M \in \mathbb{N}_0$, define the retained index set

$$\Gamma_{N,M} := \{\alpha \in \mathcal{J} : \mathrm{supp}(\alpha) \subset \{1,\dots,N\},\ |\alpha| \le M\} \tag{14}$$

The associated retained chaos space is

$$\mathcal{V}_{N,M} := \mathrm{span}\{H_\alpha(\xi) : \alpha \in \Gamma_{N,M}\}.$$

The exact retained chaos truncation is

$$U_{N,M}(t,x) := \sum_{\alpha\in\Gamma_{N,M}} c_\alpha(t,x) H_\alpha(\xi_1,\dots,\xi_N). \tag{15}$$

This notation will be used throughout the remainder of the paper.

The corresponding learned reconstruction is

$$\widehat{U}_{N,M}(t,x) := \sum_{\alpha\in\Gamma_{N,M}} \widehat{c}_\alpha(t,x)\, H_\alpha(\xi_1,\dots,\xi_N). \tag{16}$$

This representation indicates a neat reuse of the stochastic components across solutions of different SPDEs:

**Proposition 3.9** (Reusable chaos basis transfer principle)**.** *Fix an isonormal Gaussian process $W$, an orthonormal noise basis $(e_k)_{k\ge1}$, and truncation levels $N \in \mathbb{N}$ and $M \in \mathbb{N}_0$. Let $\{U^\lambda(t,x) : \lambda \in \Lambda\}$ be any family of square integrable SPDE solution fields that are measurable with respect to the same Gaussian noise. Then, for every $\lambda$, the orthogonal projection onto the retained chaos space has the representation*

$$U^\lambda_{N,M}(t,x) = \sum_{\alpha\in\Gamma_{N,M}} c^\lambda_\alpha(t,x) H_\alpha(\xi_1,\dots,\xi_N), \qquad \xi_k = W(e_k).$$

*The random variables $\{H_\alpha(\xi) : \alpha \in \Gamma_{N,M}\}$ do not depend on $\lambda$; only the deterministic coefficient fields $c^\lambda_\alpha(t,x)$ depend on the equation. Consequently, for fixed Gaussian samples $\xi^{(r)}$, the matrix*

$$\mathsf{H}_{r,\alpha} := H_\alpha(\xi^{(r)})$$

*may be computed once and reused to generate finite chaos samples for every SPDE in the family. Thus, once $W$, $(e_k)$, $N$, and $M$ are fixed, the stochastic basis is universal with respect to the family of SPDEs sharing that Gaussian driving noise.*

*Proof.* Fix $(t,x,\lambda)$ and suppress these variables temporarily. For $\alpha,\beta\in\Gamma_{N,M}$, the global definitions of $H_\alpha$ and $\alpha!$ involve only the first $N$ coordinates. Since $\xi_k=W(e_k)$ are independent standard Gaussian random variables, the tensor product Hermite polynomials satisfy

$$\mathbb{E}\big[H_\alpha(\xi)H_\beta(\xi)\big]=\prod_{k=1}^{N}\mathbb{E}\big[H_{\alpha_k}(\xi_k)H_{\beta_k}(\xi_k)\big]=\alpha!\,\mathbb{1}_{\{\alpha=\beta\}}.$$

Therefore, $\mathcal{V}_{N,M}$ is a finite dimensional subspace of $L^2(\Omega)$, spanned by mutually orthogonal Hermite basis elements.

By orthogonality of the Hermite basis, the coefficient of the orthogonal projection of $U^\lambda(t,x)$ onto $\mathcal{V}_{N,M}$ in the direction $H_\alpha(\xi)$ is the inner product with that basis element divided by its squared norm. Thus

$$c_\alpha^\lambda(t,x)=\frac{\mathbb{E}\big[U^\lambda(t,x)H_\alpha(\xi)\big]}{\mathbb{E}\big[H_\alpha(\xi)^2\big]}=\frac{1}{\alpha!}\mathbb{E}\big[U^\lambda(t,x)H_\alpha(\xi)\big],$$

and the candidate projection is

$$U_{N,M}^\lambda(t,x)=\sum_{\alpha\in\Gamma_{N,M}}\frac{\mathbb{E}[U^\lambda(t,x)H_\alpha(\xi)]}{\alpha!}H_\alpha(\xi).$$

To verify the projection property explicitly, take any $\beta\in\Gamma_{N,M}$. Orthogonality gives

$$\begin{aligned}&\mathbb{E}\left[\left(U^\lambda(t,x)-U_{N,M}^\lambda(t,x)\right)H_\beta(\xi)\right]\\&\quad=\mathbb{E}\big[U^\lambda(t,x)H_\beta(\xi)\big]-\sum_{\alpha\in\Gamma_{N,M}}\frac{\mathbb{E}[U^\lambda(t,x)H_\alpha(\xi)]}{\alpha!}\mathbb{E}[H_\alpha(\xi)H_\beta(\xi)]=0.\end{aligned}$$

Hence the residual is orthogonal to every element of $\mathcal{V}_{N,M}$, which characterizes $U_{N,M}^\lambda(t,x)$ as the orthogonal projection.

It remains to identify what is universal and what depends on the equation. The coordinates $\xi_k=W(e_k)$, and therefore all the random variables $H_\alpha(\xi)$, are fixed once $W$, $(e_k)$, $N$, and $M$ are chosen; none of them contains the equation parameter $\lambda$. All equation dependence is confined to the deterministic inner products

$$\mathbb{E}\big[U^\lambda(t,x)H_\alpha(\xi)\big]$$

and hence to the coefficient fields $c_\alpha^\lambda(t,x)$.

Finally, choose Gaussian vectors $\xi^{(1)},\dots,\xi^{(R)}$ once and set $\mathsf{H}_{r,\alpha}=H_\alpha(\xi^{(r)})$. If $\mathbf{c}^\lambda(t,x)$ is the vector of retained coefficients, then the corresponding reconstructed samples satisfy

$$\begin{pmatrix}U_{N,M}^{\lambda,(1)}(t,x)\\ \vdots\\ U_{N,M}^{\lambda,(R)}(t,x)\end{pmatrix}=\mathsf{H}\,\mathbf{c}^\lambda(t,x).$$

Changing $\lambda$ replaces only the coefficient vector; the random matrix $\mathsf{H}$ remains unchanged. This proves both the representation and the reuse assertion, and establishes the stated universality within the fixed noise coordinates and truncation levels. $\square$

This principle is conceptually simple but important: changing the SPDE transfers only the deterministic coefficient map, not the random chaos basis. The stochastic basis is universal within a family sharing the fixed Gaussian noise coordinates and truncation levels. This makes the decoupling between learning and stochastic simulation mathematically explicit.

Thus NeSTR leads to a concrete approximation procedure:

(1) choose the retained noise modes $e_1, \dots, e_N$;

(2) restrict the $S$-Transform to $h^{(N)}(a) = \sum_{k=1}^{N} a_k e_k$;

(3) learn the deterministic map $u^{(N)}(t, x; a)$;

(4) differentiate the learned map at $a = 0$ to obtain $\widehat{c}_\alpha(t, x)$;

(5) reconstruct the stochastic solution through a truncated Hermite series.

The two truncation parameters have different roles. The parameter $N$ controls how many noise directions are retained. The parameter $M$ controls how many Wiener chaos orders are kept in the reconstruction. The approximation $\widehat{U}_{N,M}$ is therefore not a single black box output: it is a structured stochastic object with explicit Hermite coefficients. This structure allows one to inspect the deterministic component, the first order response to each noise mode, and the higher order interactions separately.

**Consistency of the finite dimensional reconstruction.** We now formalize the previous discussion. The result below is deliberately stated at fixed $N$ and fixed chaos order $M$, because these are the quantities that are actually used in the computation. It says that if the learned transform has the correct derivatives at the origin, then the reconstructed stochastic field has the correct finite dimensional chaos approximation.

This statement separates two issues. The first issue is analytic: the exact $S$-Transform determines the exact chaos coefficients through Theorem 3.8. The second issue is numerical: the learned map must approximate these derivatives well enough. The next theorem addresses the second issue at fixed truncation levels, while the subsequent remark explains how this local statement fits into a full convergence picture.

**Theorem 3.10** (Consistency of the NeSTR reconstruction). *Fix $(t, x) \in [0, T] \times D$, $N \in \mathbb{N}$, and $M \in \mathbb{N}_0$, and suppose $U(t, x) \in L^2(\Omega)$. Let $u^{(N)}(t, x; a)$ be the restricted $S$-Transform (8). Suppose there exists a sequence of network parameters $(\theta_k)_{k \geq 1}$ such that*

$$\eta_k := \max_{\alpha \in \Gamma_{N,M}} \frac{1}{\alpha!} \left| \partial_a^\alpha u_{\theta_k}(t, x; 0) - \partial_a^\alpha u^{(N)}(t, x; 0) \right| \longrightarrow 0.$$

*Let $\widehat{U}_{N,M}^{(k)}(t, x)$ denote the reconstruction obtained from $u_{\theta_k}$. Then*

$$\left\| \widehat{U}_{N,M}^{(k)}(t, x) - U_{N,M}(t, x) \right\|_{L^2(\Omega)}^2 \leq \eta_k^2 \sum_{\alpha \in \Gamma_{N,M}} \alpha!,$$

*and hence $\widehat{U}_{N,M}^{(k)}(t, x) \to U_{N,M}(t, x)$ in $L^2(\Omega)$. Moreover, for any prescribed squared error tolerance $\varepsilon > 0$, it is sufficient to choose $k$ so that*

$$\eta_k \leq \left( \frac{\varepsilon}{\sum_{\alpha \in \Gamma_{N,M}} \alpha!} \right)^{1/2}.$$

*Proof.* Fix $(t, x)$ and suppress it from the notation. By Theorem 3.8, the exact coefficients are

$$c_\alpha = \frac{1}{\alpha!} \partial_a^\alpha u^{(N)}(0), \qquad \alpha \in \Gamma_{N,M}.$$

The reconstruction associated with $\theta_k$ has coefficients

$$\widehat{c}_{\alpha,k} = \frac{1}{\alpha!}\partial_a^\alpha u_{\theta_k}(0), \qquad \alpha \in \Gamma_{N,M},$$

so that $|\widehat{c}_{\alpha,k} - c_\alpha| \le \eta_k$. The Hermite polynomials are orthogonal and satisfy $\mathbb{E}[H_\alpha H_\beta] = \alpha!\,\delta_{\alpha\beta}$. Hence

$$\left\|\widehat{U}^{(k)}_{N,M} - U_{N,M}\right\|^2_{L^2(\Omega)} = \sum_{\alpha\in\Gamma_{N,M}} \alpha!\,|\widehat{c}_{\alpha,k} - c_\alpha|^2 \le \eta_k^2 \sum_{\alpha\in\Gamma_{N,M}} \alpha!.$$

Since the sum is finite and $\eta_k \to 0$, the right hand side tends to zero. Solving the same estimate for a prescribed tolerance $\varepsilon$ gives the stated epsilon condition. □

The hypothesis in the above result is explicitly about the derivatives of a sequence of trained networks. A standard universal approximation theorem on a real compact set gives uniform approximation of function values, but does not by itself imply convergence of derivatives. One possible remedy is uniform convergence of holomorphic approximants on a fixed complex polydisc around the origin; Proposition 3.12 then gives the required derivative convergence by Cauchy estimates. In the numerical scheme below, the tensor polynomial representation makes the relevant derivatives explicit and allows their errors to be checked directly.

The preceding theorem is local in the truncation parameters. The next result records the corresponding global error decomposition. It separates the part of the error caused by discarding noise modes and higher chaos orders from the part caused by learning inaccurate coefficients.

**Theorem 3.11** (Convergence with respect to mode, chaos, and learning errors)**.** *Fix $(t,x)$ and suppose*

$$U(t,x) = \sum_\alpha c_\alpha(t,x) H_\alpha(\xi) \qquad \text{in } L^2(\Omega).$$

*For $N \in \mathbb{N}$ and $M \in \mathbb{N}_0$, set $\Gamma_{N,M}$ as in (14). Let*

$$\widehat{U}_{N,M}(t,x) = \sum_{\alpha\in\Gamma_{N,M}} \widehat{c}_\alpha(t,x) H_\alpha(\xi).$$

*Then*

$$\left\|U(t,x) - \widehat{U}_{N,M}(t,x)\right\|^2_{L^2(\Omega)} = \sum_{\alpha\notin\Gamma_{N,M}} \alpha!\,|c_\alpha(t,x)|^2 + \sum_{\alpha\in\Gamma_{N,M}} \alpha!\,|c_\alpha(t,x) - \widehat{c}_\alpha(t,x)|^2. \tag{17}$$

*Consequently, if $N_j, M_j \to \infty$ and*

$$\sum_{\alpha\in\Gamma_{N_j,M_j}} \alpha!\,|c_\alpha(t,x) - \widehat{c}^{(j)}_\alpha(t,x)|^2 \longrightarrow 0,$$

*then $\widehat{U}_{N_j,M_j}(t,x) \to U(t,x)$ in $L^2(\Omega)$.*

*Proof.* Fix $(t,x)$ and suppress it from the notation. Since

$$U = \sum_\alpha c_\alpha H_\alpha \qquad \text{in } L^2(\Omega),$$

the exact projection of $U$ onto the finite set of chaos indices $\Gamma_{N,M}$ is

$$U_{N,M} := \sum_{\alpha\in\Gamma_{N,M}} c_\alpha H_\alpha.$$

The learned approximation has the same finite Hermite support, but with learned coefficients:

$$\widehat{U}_{N,M} := \sum_{\alpha \in \Gamma_{N,M}} \widehat{c}_\alpha H_\alpha.$$

We first split the error into the part outside the retained index set and the part inside it:

$$U - \widehat{U}_{N,M} = \underbrace{\sum_{\alpha \notin \Gamma_{N,M}} c_\alpha H_\alpha}_{\text{mode and chaos truncation error}} + \underbrace{\sum_{\alpha \in \Gamma_{N,M}} (c_\alpha - \widehat{c}_\alpha) H_\alpha}_{\text{learned coefficient error}}.$$

These two sums are orthogonal in $L^2(\Omega)$, because they involve disjoint sets of Hermite basis elements. More generally, for any two finite or square summable chaos series,

$$\mathbb{E}[H_\alpha H_\beta] = \alpha! \delta_{\alpha\beta}$$

implies that cross terms vanish whenever $\alpha \neq \beta$. Hence Pythagoras' identity gives

$$\left\| U - \widehat{U}_{N,M} \right\|_{L^2(\Omega)}^2 = \| U - U_{N,M} \|_{L^2(\Omega)}^2 + \left\| U_{N,M} - \widehat{U}_{N,M} \right\|_{L^2(\Omega)}^2.$$

Applying the Hermite norm identity separately to the two terms yields

$$\| U - U_{N,M} \|_{L^2(\Omega)}^2 = \sum_{\alpha \notin \Gamma_{N,M}} \alpha! \, |c_\alpha|^2$$

and

$$\left\| U_{N,M} - \widehat{U}_{N,M} \right\|_{L^2(\Omega)}^2 = \sum_{\alpha \in \Gamma_{N,M}} \alpha! \, |c_\alpha - \widehat{c}_\alpha|^2.$$

Combining these two identities proves (17).

It remains to justify the convergence statement. Because $U \in L^2(\Omega)$, its chaos energy is finite:

$$\sum_\alpha \alpha! \, |c_\alpha|^2 < \infty.$$

The sets $\Gamma_{N_j,M_j}$ exhaust all finite multi-indices as $N_j, M_j \to \infty$: for any fixed multi-index $\alpha$, there is an index $j_0$ such that $\alpha \in \Gamma_{N_j,M_j}$ for all $j \geq j_0$. Therefore the tails

$$\sum_{\alpha \notin \Gamma_{N_j,M_j}} \alpha! \, |c_\alpha|^2$$

tend to zero by monotone convergence for series.The remaining term in (17) is exactly the learned coefficient error over the retained set, and this term tends to zero by assumption. Thus both contributions to the squared $L^2(\Omega)$ error vanish, and consequently

$$\left\| U - \widehat{U}_{N_j,M_j} \right\|_{L^2(\Omega)} \longrightarrow 0.$$

This proves the claimed convergence. □

The coefficient error in Theorem 3.11 can be controlled directly from a uniform transform error on a small complex neighborhood of the origin. This is useful because it explains why sampling near $a = 0$ matters.

**Proposition 3.12** (Derivative stability from local transform accuracy)**.** *Fix* $(t,x)$*,* $N \in \mathbb{N}$*, and* $M \in \mathbb{N}_0$*. Let* $P_r := \{a \in \mathbb{C}^N : |a_k| \le r,\ k = 1,\dots,N\}$ *be a polydisc centered at zero. Assume that* $u^{(N)}(t,x;\cdot)$ *and* $u_\theta(t,x;\cdot)$ *are holomorphic on a neighborhood of* $P_r$ *and that*

$$\sup_{a\in P_r} |u_\theta(t,x;a) - u^{(N)}(t,x;a)| \le \varepsilon.$$

*Then, for every* $\alpha \in \mathbb{N}_0^N$ *with* $|\alpha| \le M$*,*

$$|\widehat{c}_\alpha(t,x) - c_\alpha(t,x)| \le \frac{\varepsilon}{r^{|\alpha|}}. \tag{18}$$

*In particular,*

$$\left\|\widehat{U}_{N,M}(t,x) - U_{N,M}(t,x)\right\|^2_{L^2(\Omega)} \le \varepsilon^2 \sum_{\alpha\in\Gamma_{N,M}} \frac{\alpha!}{r^{2|\alpha|}}. \tag{19}$$

*Proof.* Fix $(t,x)$ and suppress it from the notation. Define the transform error

$$e(a) := u_\theta(a) - u^{(N)}(a).$$

By assumption, $e$ is holomorphic on a neighborhood of the closed polydisc $P_r$, and

$$\sup_{a\in P_r} |e(a)| \le \varepsilon.$$

Let $\alpha = (\alpha_1,\dots,\alpha_N)$. For any radius $0 < \rho < r$, the polydisc Cauchy formula gives

$$\partial_a^\alpha e(0) = \frac{\alpha!}{(2\pi i)^N} \int_{|z_1|=\rho} \cdots \int_{|z_N|=\rho} \frac{e(z_1,\dots,z_N)}{z_1^{\alpha_1+1}\cdots z_N^{\alpha_N+1}}\, dz_N \cdots dz_1.$$

Taking absolute values and using $|z_k| = \rho$ on each contour yields

$$\begin{aligned}
|\partial_a^\alpha e(0)| &\le \frac{\alpha!}{(2\pi)^N}\varepsilon \prod_{k=1}^N \int_{|z_k|=\rho} \frac{|dz_k|}{|z_k|^{\alpha_k+1}} \\
&= \frac{\alpha!}{(2\pi)^N}\varepsilon \prod_{k=1}^N \frac{2\pi\rho}{\rho^{\alpha_k+1}} = \alpha!\,\varepsilon\,\rho^{-|\alpha|}.
\end{aligned}$$

Since this estimate holds for every $\rho < r$, letting $\rho \uparrow r$ gives

$$|\partial_a^\alpha e(0)| \le \alpha!\,\varepsilon\, r^{-|\alpha|}.$$

The recovered exact and learned chaos coefficients are

$$c_\alpha = \frac{1}{\alpha!}\partial_a^\alpha u^{(N)}(0), \qquad \widehat{c}_\alpha = \frac{1}{\alpha!}\partial_a^\alpha u_\theta(0).$$

Therefore

$$\widehat{c}_\alpha - c_\alpha = \frac{1}{\alpha!}\partial_a^\alpha e(0),$$

and the derivative estimate gives

$$|\widehat{c}_\alpha - c_\alpha| \le \varepsilon r^{-|\alpha|},$$

which proves (18).

It remains to translate this coefficient estimate into a reconstruction estimate. Using the exact retained chaos truncation $U_{N,M}$ from (15), orthogonality of the Hermite basis gives

$$\left\|\widehat{U}_{N,M} - U_{N,M}\right\|_{L^2(\Omega)}^2 = \sum_{\alpha\in\Gamma_{N,M}} \alpha!\,|\widehat{c}_\alpha - c_\alpha|^2.$$

Substituting the coefficient bound just proved into this finite sum gives

$$\left\|\widehat{U}_{N,M} - U_{N,M}\right\|_{L^2(\Omega)}^2 \le \sum_{\alpha\in\Gamma_{N,M}} \alpha!\,\varepsilon^2 r^{-2|\alpha|} = \varepsilon^2 \sum_{\alpha\in\Gamma_{N,M}} \frac{\alpha!}{r^{2|\alpha|}},$$

which is (19). $\square$

The next result wraps up convergence bounds across chaos truncations, numerical solver and neural network approximation:

**Theorem 3.13** (NeSTR convergence bound)**.** *Fix $(t,x)$, truncation levels $N\in\mathbb{N}$ and $M\in\mathbb{N}_0$, and a radius $r>0$. Let*

$$U_N(t,x) := \sum_{\substack{\alpha\in\mathcal{J}\\ \mathrm{supp}(\alpha)\subset\{1,\dots,N\}}} c_\alpha(t,x) H_\alpha(\xi)$$

*be the full chaos projection onto the first $N$ noise modes, and let $U_{N,M}(t,x)$ be its order-$M$ truncation from (15). Define*

$$\varepsilon_{\text{mode}} := \|U(t,x) - U_N(t,x)\|_{L^2(\Omega)}, \qquad \varepsilon_{\text{chaos}} := \|U_N(t,x) - U_{N,M}(t,x)\|_{L^2(\Omega)}.$$

*Let $\widetilde{u}^{(N)}$ be a deterministic numerical approximation of the exact restricted transform $u^{(N)}$, and let $u_\theta$ be the trained NeSTR approximation. Assume that these maps are holomorphic on a neighborhood of the polydisc $P_r$ and satisfy*

$$\sup_{a\in P_r} |\widetilde{u}^{(N)}(a) - u^{(N)}(a)| \le \varepsilon_{\text{solver}}, \qquad \sup_{a\in P_r} |u_\theta(a) - \widetilde{u}^{(N)}(a)| \le \varepsilon_{\text{NN}}.$$

*Then the reconstructed field satisfies*

$$\left\|U - \widehat{U}_{N,M}\right\|_{L^2(\Omega)} \le \varepsilon_{\text{mode}} + \varepsilon_{\text{chaos}} + C_{N,M}(r)(\varepsilon_{\text{solver}} + \varepsilon_{\text{NN}}), \tag{20}$$

*where*

$$C_{N,M}(r) := \left(\sum_{\alpha\in\Gamma_{N,M}} \frac{\alpha!}{r^{2|\alpha|}}\right)^{1/2}.$$

*Consequently, along any sequence for which the mode and chaos errors vanish and $C_{N,M}(r)(\varepsilon_{\text{solver}}+\varepsilon_{\text{NN}})\to 0$, the NeSTR reconstruction converges to $U$ in $L^2(\Omega)$.*

*Proof.* Fix $(t,x)$ and suppress it from the notation. We first separate the two deterministic transform errors before returning to the stochastic field. Define

$$e_{\text{solver}}(a) := \widetilde{u}^{(N)}(a) - u^{(N)}(a), \qquad e_{\text{NN}}(a) := u_\theta(a) - \widetilde{u}^{(N)}(a).$$

Both functions are holomorphic on a neighborhood of $P_r$, and the assumptions give

$$\sup_{a\in P_r} |e_{\text{solver}}(a)| \le \varepsilon_{\text{solver}}, \qquad \sup_{a\in P_r} |e_{\text{NN}}(a)| \le \varepsilon_{\text{NN}}.$$

Introduce the intermediate coefficients recovered from the deterministic solver,

$$\widetilde{c}_\alpha := \frac{1}{\alpha!}\partial_a^\alpha \widetilde{u}^{(N)}(0).$$

The exact and trained coefficients are, respectively,

$$c_\alpha = \frac{1}{\alpha!}\partial_a^\alpha u^{(N)}(0), \qquad \widehat{c}_\alpha = \frac{1}{\alpha!}\partial_a^\alpha u_\theta(0).$$

Applying the coefficient estimate (18) separately to $e_{\text{solver}}$ and $e_{\text{NN}}$ yields, for every $\alpha \in \Gamma_{N,M}$,

$$|\widetilde{c}_\alpha - c_\alpha| \le \varepsilon_{\text{solver}} r^{-|\alpha|}, \qquad |\widehat{c}_\alpha - \widetilde{c}_\alpha| \le \varepsilon_{\text{NN}} r^{-|\alpha|}.$$

The triangle inequality at the coefficient level therefore gives

$$|\widehat{c}_\alpha - c_\alpha| \le \big(\varepsilon_{\text{solver}} + \varepsilon_{\text{NN}}\big) r^{-|\alpha|}. \tag{21}$$

This step shows explicitly how deterministic solver error and neural learning error enter the same recovered coefficient without being identified with one another.

We next pass from coefficient errors to the retained stochastic field. From the definitions of $U_{N,M}$ and $\widehat{U}_{N,M}$,

$$U_{N,M} - \widehat{U}_{N,M} = \sum_{\alpha \in \Gamma_{N,M}} (c_\alpha - \widehat{c}_\alpha) H_\alpha(\xi).$$

Hermite orthogonality eliminates all cross terms. Consequently,

$$\begin{aligned}\left\| U_{N,M} - \widehat{U}_{N,M} \right\|_{L^2(\Omega)}^2 &= \sum_{\alpha \in \Gamma_{N,M}} \alpha!\, |c_\alpha - \widehat{c}_\alpha|^2 \\ &\le \big(\varepsilon_{\text{solver}} + \varepsilon_{\text{NN}}\big)^2 \sum_{\alpha \in \Gamma_{N,M}} \frac{\alpha!}{r^{2|\alpha|}}.\end{aligned}$$

Taking square roots and using the definition of $C_{N,M}(r)$ gives

$$\left\| U_{N,M} - \widehat{U}_{N,M} \right\|_{L^2(\Omega)} \le C_{N,M}(r)\big(\varepsilon_{\text{solver}} + \varepsilon_{\text{NN}}\big). \tag{22}$$

For completeness, the three stochastic errors lie in mutually orthogonal chaos subspaces. Indeed, $U_N$ is the orthogonal projection onto the full chaos space generated by $\xi_1, \ldots, \xi_N$, so $U - U_N$ is orthogonal to that space. The difference $U_N - U_{N,M}$ contains only chaos components in the first $N$ modes whose order is larger than $M$. In contrast, $U_{N,M} - \widehat{U}_{N,M}$ belongs to the retained space $\mathcal{V}_{N,M}$. Pythagoras' identity therefore gives the stronger intermediate estimate

$$\begin{aligned}\left\| U - \widehat{U}_{N,M} \right\|_{L^2(\Omega)}^2 &= \|U - U_N\|_{L^2(\Omega)}^2 + \|U_N - U_{N,M}\|_{L^2(\Omega)}^2 \\ &\quad + \left\| U_{N,M} - \widehat{U}_{N,M} \right\|_{L^2(\Omega)}^2 \\ &\le \varepsilon_{\text{mode}}^2 + \varepsilon_{\text{chaos}}^2 + C_{N,M}(r)^2 \big(\varepsilon_{\text{solver}} + \varepsilon_{\text{NN}}\big)^2.\end{aligned}$$

Taking square roots and using $\sqrt{x^2 + y^2 + z^2} \le x + y + z$ for nonnegative $x, y, z$ proves (20).

Finally, along the sequence in the statement, $\varepsilon_{\text{mode}} \to 0$, $\varepsilon_{\text{chaos}} \to 0$, and

$$C_{N,M}(r)\big(\varepsilon_{\text{solver}} + \varepsilon_{\text{NN}}\big) \longrightarrow 0.$$

Every term on the right hand side of (20) tends to zero. Hence $\widehat{U}_{N,M} \to U$ in $L^2(\Omega)$, completing the proof. □

**Corollary 3.14** (Convergence of NeSTR reconstructions)**.** *Assume that a sequence of learned transforms $u_{\theta_n}$ converges to $u^{(N)}$ in a topology strong enough to imply convergence of the derivatives $\partial_a^\alpha u_{\theta_n}(t, x; 0)$ for all $|\alpha| \le M$. Then the corresponding reconstructions $\widehat{U}^{(n)}_{N,M}(t, x)$ converge in $L^2(\Omega)$ to the exact $N$-mode, $M$-th order chaos truncation. If, in addition, there exist U-functionals $\Phi_n, \Phi$ on $E_{\mathbb{C}}$ whose restrictions satisfy*

$$\Phi_n(h^{(N)}(a)) = u_{\theta_n}(t, x; a), \qquad \Phi(h^{(N)}(a)) = u^{(N)}(t, x; a),$$

*and $\Phi_n \to \Phi$ in the Potthoff–Streit U-functional topology of Theorem 2.2, then the corresponding Hida distributions $F_n = S^{-1}\Phi_n$ converge to $F = S^{-1}\Phi$ in $(\mathcal{S})^*$.*

*Proof.* Fix $(t, x)$ and suppress this dependence in the notation. Recall the exact retained chaos truncation $U_{N,M}$ from (15). By the inverse $S$-Transform identity,

$$c_\alpha = \frac{1}{\alpha!}\partial_a^\alpha u^{(N)}(0), \qquad |\alpha| \le M.$$

The reconstruction obtained from the learned transform has coefficients

$$\widehat{c}_{\alpha,n} = \frac{1}{\alpha!}\partial_a^\alpha u_{\theta_n}(0), \qquad |\alpha| \le M.$$

The assumed convergence of the learned transforms in a topology controlling these derivatives therefore implies

$$\widehat{c}_{\alpha,n} \longrightarrow c_\alpha \qquad \text{for every } \alpha \in \mathbb{N}_0^N \text{ with } |\alpha| \le M.$$

Since there are only finitely many such multi-indices, coefficient-wise convergence is enough to give convergence of the finite Hermite expansion. More explicitly, using the orthogonality relation $\mathbb{E}[H_\alpha H_\beta] = \alpha!\delta_{\alpha\beta}$, we obtain

$$\left\|\widehat{U}^{(n)}_{N,M} - U_{N,M}\right\|^2_{L^2(\Omega)} = \mathbb{E}\left[\left|\sum_{\substack{\alpha\in\mathbb{N}_0^N\\|\alpha|\le M}} (\widehat{c}_{\alpha,n} - c_\alpha)H_\alpha(\xi_1, \dots, \xi_N)\right|^2\right] = \sum_{\substack{\alpha\in\mathbb{N}_0^N\\|\alpha|\le M}} \alpha!\,|\widehat{c}_{\alpha,n} - c_\alpha|^2.$$

Each term in the last finite sum tends to zero, hence the whole sum tends to zero. This proves convergence in $L^2(\Omega)$ to the exact finite dimensional, finite order chaos truncation.

We now prove the assertion in the Hida distribution topology under the stronger U-functional hypothesis. A restricted map on $\mathbb{C}^N$ is not, by itself, a full $S$-Transform on $E_{\mathbb{C}}$; the assumed extensions $\Phi_n$ and $\Phi$ are therefore essential. By Theorem 2.2, define

$$F_n := S^{-1}\Phi_n, \qquad F := S^{-1}\Phi.$$

The hypothesis is precisely

$$\Phi_n \longrightarrow \Phi$$

in the Potthoff–Streit U-functional topology of Theorem 2.2. By that theorem, the $S$-Transform is a topological isomorphism from $(\mathcal{S})^*$ onto the U-functional space. Applying the continuous inverse map gives

$$F_n = S^{-1}(\Phi_n) \longrightarrow S^{-1}(\Phi) = F \qquad \text{in } (\mathcal{S})^*.$$

Thus convergence of deterministic transforms in the U-functional topology transfers back to convergence of the corresponding stochastic objects in the Hida distribution topology. The two parts together show that the computable finite chaos reconstructions converge in $L^2$, while the untruncated transform objects converge in $(\mathcal{S})^*$ only under this stronger, explicitly stated U-functional convergence assumption. □

*Remark* 3.15. Theorem 2.2 is the functional analytic justification for the NeSTR learning step: the $S$-Transform is a continuous embedding and its inverse is continuous on the image. Thus, approximation of the deterministic $S$-Transform in the correct transform topology implies approximation of the stochastic object. Theorem 3.10 is the computable, truncated version of this statement. In a full convergence theorem one must additionally control:

(i) the basis truncation error $U - U_N$,

(ii) the chaos truncation error $U_N - U_{N,M}$,

(iii) the deterministic solver error $u^{(N)} - \widetilde{u}^{(N)}$,

(iv) the neural approximation error $\widetilde{u}^{(N)} - u_\theta$.

Theorem 3.13 combines these contributions quantitatively.

**Approximation theoretic perspective.** The previous constructions reduce the original SPDE problem to the approximation of the deterministic family

$$(t, x, a) \longmapsto u^{(N)}(t, x; a).$$

From the point of view of approximation theory, there are two distinct layers:

(1) Transform space restriction. The infinite dimensional transform variable $h \in \mathcal{H}_\mathbb{C}$ at the square integrable level, or $h \in E_\mathbb{C}$ in the Hida setting, is replaced by the finite dimensional coordinate vector $a \in \mathbb{C}^N$.

(2) Function approximation. The restricted transform $u^{(N)}$ is approximated by a neural network $u_\theta$ on compact subsets of $[0, T] \times D \times \mathbb{C}^N$.

If one works on a compact set

$$K \subset [0, T] \times D \times \mathbb{C}^N,$$

identifying $\mathbb{C}^N$ with $\mathbb{R}^{2N}$, standard universal approximation theorems on finite dimensional Euclidean domains [8, 16] imply that continuous neural networks can approximate $u^{(N)}$ uniformly on $K$, provided the architecture is sufficiently rich. For the additive noise heat equation this is immediate, since $u^{(N)}$ is affine linear in $a$. For nonlinear transformed equations, one expects the same conclusion once continuity of $u^{(N)}$ with respect to $(t, x, a)$ has been established.

*Remark* 3.16. This finite dimensional approximation step is the point at which ideas from operator learning become relevant. The transformed SPDE induces an operator from a finite dimensional representation of the noise parameter $a$ to a deterministic solution field. In this sense, NeSTR learns a restricted solution operator in transform space.

# 4 NeSTR numerical results

The numerical section is organized as a sequence of increasingly demanding tests. The first example is deterministic and checks the underlying backward solver in the absence of noise using the benchmark problem (24). The second example is the additive noise heat equation (25), where the restricted transform is affine in $a$ and the chaos coefficients are known in closed form (Proposition 3.5). This validates the full NeSTR pipeline in a controlled setting: restriction, learning in transform space, differentiation at the origin, and chaos reconstruction. The third example is the multiplicative Wick heat equation (26), where $a \mapsto u^{(N)}$ is genuinely nonlinear and higher chaos coefficients no longer vanish. It is therefore the decisive numerical test for high order inverse reconstruction.

### 4.1 Numerical approximation scheme and implementation details

This subsection describes the implementation used in the experiments. The focus is on the backward dynamic programming solver, the polynomial representation in transform variables, and the way coefficient derivatives are extracted.

After the $S$-Transform, the equation for $u^{(N)}(t,x;a)$ is a deterministic parabolic Cauchy problem in $(t,x)$ that carries the finite dimensional transform coordinate $a$ as a fixed parameter; in the additive case it is (4). We do not discretize the physical space on a grid. A mesh based solver would have to resolve the transform coordinates as well as the physical ones, since $u^{(N)}$ must be known as a function of $a$ in order to be differentiated at the origin, and a grid over $[-A,A]^N$ is not affordable already for moderate $N$. Instead we use its Feynman–Kac representation (6) and approximate $u_\theta$ by a deep backward dynamic programming scheme for the associated backward stochastic differential equation, following the deep PDE/BSDE learning philosophy of [12, 3, 2]: the time interval is split into steps and one network is trained per time level, backwards from the terminal datum. Since $a$ enters only the driver and not the driving noise, the simulated diffusion stays in the physical spatial dimension $d$. Increasing $N$ therefore does not enlarge the simulated diffusion state, although it does increase the parametric approximation cost. Among the probabilistic schemes, backward dynamic programming is the one that returns what the recovery step needs. It produces $u_\theta(t_m,\cdot;\cdot)$ as a network at every time level, so the map $a\mapsto u_\theta(t_m,x;a)$ is available in closed form and can be differentiated at $a=0$ at any $(t_m,x)$. A deep BSDE scheme in the form of [12] instead returns the value at a single initial point together with the gradient process along the sampled paths, and would have to be re-solved for every point at which a coefficient field is wanted. A physics-informed formulation would be a further alternative, but it reaches the coefficients by automatic differentiation of a residual in $a$, which is the step we avoid for the reasons given below. We emphasize that this is a choice of implementation and not of framework: as noted in Remark 3.7, any deterministic solver for the transformed equation may be used.

The stochastic benchmark equations are written as initial value problems, whereas the Deep backward dynamic programming (DBDP) solver is formulated backward from terminal data. We reconcile the two conventions by introducing the reverse time field

$$v(s,x;a) := u^{(N)}(T-s,x;a).$$

Then $v(T,x;a)=u^{(N)}(0,x;a)$, so the known initial condition becomes the terminal condition for the backward solver. After training, results are reported in the original physical time variable $t=T-s$. Below, the time grid in the algorithm refers to the backward solver coordinate, while all figures and tables use physical time.

The one nonstandard choice concerns the recovery step (13). Rather than differentiating a generic network in $a$ by automatic differentiation, which is delicate at high order, we take $u_\theta$ to be a polynomial in $a$: a tensor expansion in $a_1,\dots,a_N$ whose coefficients are functions of $(t,x)$ represented by an ordinary network. Two bases are used below, tensor Chebyshev [27] in the additive experiment and tensor monomial in the multiplicative one; both span the same polynomial space, so the recovery step is exact either way. The map $a\mapsto u_\theta(t,x;a)$ is then a genuine polynomial, so its derivatives at $a=0$, and hence the approximate coefficients $\widehat{c}_\alpha(t,x)$, are obtained exactly. This matches the analytic structure of $u^{(N)}$, which by Proposition 3.3 is itself entire in $a$.

Derivative reliability is essential here and is not a consequence of ordinary universal approximation of function values. Repeated automatic differentiation of a generic network can amplify small value errors, and a small training loss does not by itself guarantee accurate high order derivatives at $a=0$. NeSTR addresses this issue in two complementary ways. Analytically, Proposition 3.12 and Theorem 3.13 show that uniform transform accuracy on a complex neighborhood controls ev-ery retained derivative. Numerically, the dependence on $a$ is represented by a finite polynomial. Therefore $\partial_a^\alpha u_\theta(0)$ is evaluated as an exact linear combination of learned spectral coefficients, without finite differences and without repeated automatic differentiation through the spatial network. Scaling the transform box to $[-1,1]^N$, keeping the polynomial degree moderate, and checking coefficient errors order by order are consequently part of the stability strategy rather than optional implementation details. The deterioration observed at unresolved high chaos orders remains visible and is used to select the final truncation order.

The full NeSTR computational pipeline is summarized in Algorithm 1. The important point is that the neural network is trained only on deterministic transform data; stochastic information enters again only through the final inverse $S$-Transform step.

The symbol $N$ is reserved for the number of retained noise modes, so the number of time steps is written $N_t$, while $M$ denotes the chaos truncation order and $\Gamma_{N,M}$ is the retained chaos index set. The tensor degree set used below is $\Lambda_p := \{\nu \in \mathbb{N}_0^N : |\nu| \le p\}$, where $|\nu| = \nu_1 + \cdots + \nu_N$, so that $|\Lambda_p| = \binom{N+p}{p}$; for $\nu \in \Lambda_p$ we write $P_\nu$ for the corresponding tensor polynomial in the scaled variable $a/A$. Two choices of $P_\nu$ occur below. The additive experiment uses the tensor Chebyshev polynomials

$$T_\nu(a) := \prod_{k=1}^{N} T_{\nu_k}(a_k/A), \qquad T_j(\cos\vartheta) = \cos(j\vartheta),$$

where $T_j$ is the Chebyshev polynomial of the first kind of degree $j$; the multiplicative experiment uses the tensor monomials

$$(a/A)^\nu := \prod_{k=1}^{N} (a_k/A)^{\nu_k}.$$

Both families span the polynomials of total degree at most $p$, so $P_\nu$ below denotes either one and the recovery step (13) is exact in both cases. They differ in conditioning, and the monomial choice has the further feature that the coefficient of $(a/A)^\nu$ is the recovered chaos coefficient itself, up to the factorial normalization.

**Implementation details.** We use the same computational template in all numerical examples. The physical time interval is divided into a finite grid

$$0 = t_0 < t_1 < \cdots < t_{N_t} = T,$$

and the transformed deterministic equation is solved backward from the terminal or final condition. At each time level $t_m$, a separate neural approximation is trained and initialized from the network obtained at the neighboring later time level. This warm start is important in practice: adjacent time levels represent nearby deterministic solution maps, so reusing the previous parameters reduces the number of iterations needed after the first step.

The trainable ansatz has the form

$$u_\theta(t_m, x; a) = \sum_{\nu\in\Lambda_p} q_{\theta,\nu}^{(m)}(x)\, P_\nu(a), \tag{23}$$

where $P_\nu$ and $\Lambda_p$ are the tensor polynomial basis and degree set defined above, and each coefficient function $q_{\theta,\nu}^{(m)}$ is represented by a fully connected feedforward neural network in the physical variable $x$. The deterministic benchmark corresponds to the degenerate convention $N = 0$, with $h^{(0)} = 0$ and $\Gamma_{0,M} = \{0\}$. In the additive experiment below, $N = 3$ and $A = 1$; the same architecture is used for all retained modes. The polynomial representation is used only in the transform variables.

**Algorithm 1** NeSTR training and inverse chaos reconstruction

**Require:** SPDE in Wick form, initial data $u_0$, noise modes $e_1, \dots, e_N$, chaos order $M$, polynomial degree $p$, transform box $[-A, A]^N$, and time grid $0 = t_0 < \cdots < t_{N_t} = T$.

**Ensure:** Approximate chaos coefficients $\widehat{c}_\alpha(t_m, x)$ and stochastic reconstruction $\widehat{U}_{N,M}(t_m, x)$.

1: Apply the $S$-Transform to the SPDE and use $S(F \diamond G) = (SF)(SG)$ to obtain the deterministic parametric PDE for $u^{(N)}(t, x; a)$.

2: Restrict the test function to $h^{(N)}(a) = \sum_{k=1}^N a_k e_k$, with $a \in [-A, A]^N$.

3: **for** $m = N_t - 1, N_t - 2, \dots, 0$ **do**

4: Sample physical states from the Feynman–Kac/BSDE diffusion and sample transform coordinates $a$ from $[-A, A]^N$.

5: Train the polynomial–neural ansatz

$$u_\theta(t_m, x; a) = \sum_{\nu \in \Lambda_p} q_{\theta,\nu}^{(m)}(x) P_\nu(a)$$

by minimizing the one step dynamic programming loss.

6: Warm start the next backward level with the parameters learned at time $t_m$.

7: **end for**

8: **for** each time level $t_m$ and retained multi-index $\alpha \in \Gamma_{N,M}$ **do**

9: Recover the chaos coefficient analytically from the polynomial part,

$$\widehat{c}_\alpha(t_m, x) = \frac{1}{\alpha!} \partial_a^\alpha u_\theta(t_m, x; 0).$$

10: **end for**

11: Reconstruct the stochastic field by the truncated inverse $S$-Transform,

$$\widehat{U}_{N,M}(t_m, x) = \sum_{\alpha \in \Gamma_{N,M}} \widehat{c}_\alpha(t_m, x) H_\alpha(\xi).$$

The dependence on $x$ is still learned by the neural network, and the time dependence is handled by backward stepping.

Training points are generated by sampling the physical state from the Feynman–Kac/BSDE diffusion and drawing the transform coordinate $a$ from $[-A, A]^N$. The loss minimized at each time level is the empirical version of the one step dynamic programming error; we refer to this quantity as the DBDP loss. More precisely, if $\mathcal{T}_m^{(j)}$ denotes the one step Feynman–Kac or BSDE target computed from the already trained approximation at time $t_{m+1}$, then

$$\mathcal{L}_m^{\mathrm{DBDP}}(\theta) := \frac{1}{J}\sum_{j=1}^{J}\left|u_\theta(t_m, X_m^{(j)}; a^{(j)}) - \mathcal{T}_m^{(j)}\right|^2.$$

For a semilinear transformed equation, the target has the usual Euler form

$$\mathcal{T}_m^{(j)} = \widehat{u}_{m+1}(X_{m+1}^{(j)}; a^{(j)}) + f\Big(t_m, X_m^{(j)}, \widehat{u}_{m+1}(X_{m+1}^{(j)}; a^{(j)}), a^{(j)}\Big)\,\Delta t_m,$$

with the sign adjusted to the convention of the transformed PDE. Thus the DBDP curves reported below are training curves for this one step squared residual.

For examples with a closed form solution, we also evaluate the independent supervised error

$$\left(\frac{1}{J}\sum_{j=1}^{J}|u_\theta(t, x_j; a_j) - u^{(N)}(t, x_j; a_j)|^2\right)^{1/2}.$$

The reported coefficient errors are computed after differentiating (23) analytically at $a = 0$. Thus the inverse step is not affected by numerical finite differences in the transform coordinates.

For reproducibility, Table 1 records the hyperparameters used for the runs reported below. The same spatial network is used for all polynomial coefficient functions $q_{\theta,\nu}^{(m)}$ at a fixed time level; only the number of retained transform modes and the polynomial degree differ between the additive and multiplicative tests.

### 4.2 Numerical examples

This subsection reports the numerical tests in increasing order of stochastic complexity. The examples separate deterministic solver accuracy, additive noise reconstruction, and higher order chaos recovery in the multiplicative Wick case. The code corresponding to the implementations presented in this section is available on https://github.com/janrems/SPDE-S-T.

**Example 1: Deterministic benchmark.** As a control we first solve the backward heat equation

$$\begin{cases} \partial_t u(t,x) + \frac{1}{2}\Delta u(t,x) = 0, & (t,x) \in [0,T) \times \mathbb{R}^d, \\ u(T,x) = g(x) = |x|^2, & x \in \mathbb{R}^d. \end{cases} \tag{24}$$

The exact solution is

$$u(t,x) = |x|^2 + d\,(T - t).$$

We take $d = 1$ and start the simulated diffusion at $x_0 = 1$. There is no noise and no transform parameter $a$, so this example isolates the backward solver before any $S$-Transform or inverse reconstruction is involved. It reproduces the closed form to a relative $L^2$ error of at most 1.1%, with mean 0.8% across interior time levels, over the support of the simulated diffusion.

| quantity | deterministic | additive Wick | multiplicative Wick |
|---|---|---|---|
| spatial dimension $d$ | 1 | 1 | 1 |
| terminal time $T$ | 1 | 1 | 1 |
| time steps $N_t$ | 20 | 10 | 20 |
| training samples per step | 4096 | 512 | 512 |
| retained noise modes $N$ | 0 | 3 | 1 |
| transform box $[-A, A]^N$ | – | $[-1, 1]^3$ | $[-1, 1]$ |
| polynomial basis $P_\nu$ | – | $T_\nu$ | $(a/A)^\nu$ |
| polynomial total degree $p$ | – | 3 | 6 |
| chaos orders evaluated | 0 | 0–2 | 0–4 |
| hidden layers | 3 | 2 | 2 |
| neurons per hidden layer | 64 | 32 | 64 |
| activation | tanh | tanh | tanh |
| optimizer | Adam | Adam | Adam |
| learning rate | $10^{-3}$ | $10^{-3}$ | $10^{-3}$ |
| iterations per backward step | 3000 | 5000 | 20000 |
| warm start across time levels | yes | yes | yes |
| random seed | 1234 | 0 | 0 |

Table 1: Hyperparameters for the reported numerical experiments. The multiplicative run uses one retained transform mode and resolves chaos orders up to four; the usable reconstruction order is selected afterward from the measured coefficient errors. The learning rate is held fixed, and the transform coordinate is drawn uniformly from $[-A, A]^N$.

**Example 2: Additive noise heat equation.** The second test is the additive heat equation in one space dimension on the periodic domain $D = \mathbb{T} \cong [-\pi, \pi]$, where the endpoints are identified. Thus the boundary conditions are explicitly

$$U(t, -\pi) = U(t, \pi), \qquad \partial_x U(t, -\pi) = \partial_x U(t, \pi).$$

The equation is

$$\begin{cases} \partial_t U(t, x) = \Delta U(t, x) + \dot{W}(t, x), & (t, x) \in (0, 1] \times D, \\ U(0, x) = \sin x, & x \in D. \end{cases} \tag{25}$$

For this benchmark the isonormal Gaussian process is defined over

$$\mathcal{H} = L^2\left([0, 1] \times \mathbb{T}, ds\, \frac{dx}{2\pi}\right).$$

We retain three modes that are constant in time and have the spatial profiles

$$\phi_1(x) = \sqrt{2}\, \sin x, \qquad \phi_2(x) = \sqrt{2}\, \cos x, \qquad \phi_3(x) = \sqrt{2}\, \sin(2x),$$

with Laplacian eigenvalues $\lambda_1 = \lambda_2 = 1$ and $\lambda_3 = 4$. With the uniform torus probability measure $dx/(2\pi)$, the space–time modes $e_k(s, x) = \phi_k(x)$ are orthonormal in $\mathcal{H}$. Thus, the corresponding coordinates $\xi_k = W(e_k)$ are independent standard Gaussian variables, in agreement with the Hermite normalization used in the reconstruction.

We sample the transform parameter from $[-1, 1]^3$. The exact restricted transform has the form

$$u^{(N)}(t, x; a) = c_0(t, x) + \sum_{k=1}^{3} a_k c_k(t, x), \qquad c_0(t, x) = e^{-t} \sin x, \qquad c_k(t, x) = \phi_k(x) \frac{1 - e^{-\lambda_k t}}{\lambda_k}.$$

Thus the only nonvanishing retained coefficients are $c_0$ and the three first order coefficients $c_k$; every coefficient of order $\geq 2$ is identically zero. This gives an exact target for both NeSTR learning and inverse reconstruction. Figure 1 shows the exact restricted transform at the coordinate vector $a = (0.5, -0.5, 0.5)$.

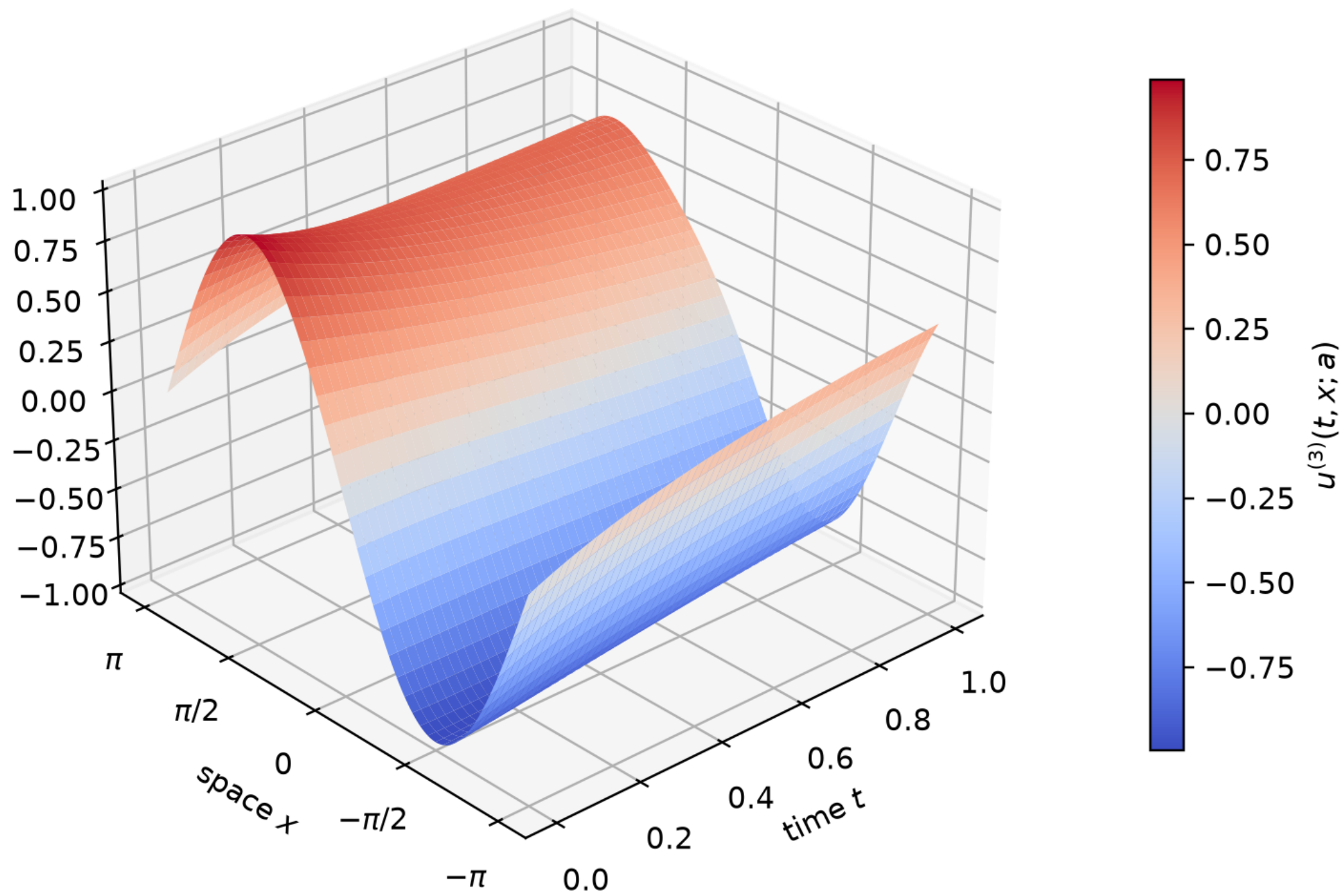

Figure 1: Exact additive restricted $S$-Transform surface at $a = (0.5, -0.5, 0.5)$ on the periodic spatial domain $D = [-\pi, \pi]$.

Figure 2 tests the learned dependence on the transform coordinates. In each panel, all coordinates except one are fixed at zero and the learned map $u_\theta(t, x; a)$ is plotted as a function of the remaining coordinate. For the additive equation the exact $S$-Transform is affine in $a$, so each slice should be a straight line. The slope at the origin is the corresponding first order chaos coefficient. The learned curves follow the closed form solution closely; the mild curvature visible near the boundary of the interval $[-1, 1]$ represents leakage into higher order chaos terms, which should vanish for this test problem.

Figure 3 compares the learned restricted transform with the closed form solution over the $(t, x)$-plane for one fixed value of $a$. The left and middle panels show that the network captures the global shape of the solution. The right panel displays the signed error, which is largest near the spatial edges; this is expected because the simulated diffusion visits those regions less frequently. Figure 4 then evaluates the inverse step. It compares the recovered coefficient fields with the exact coefficients at a fixed time. The deterministic coefficient and the three first order modes are all recovered with the correct qualitative shape, while the remaining error is mainly an amplitude error in the oscillatory modes.

Table 2 and Figure 5 summarize the accuracy across time. We report absolute root mean square errors rather than relative errors. This is important because the first order response $c_k = \phi_k(x)(1 - e^{-\lambda_k t})/\lambda_k$ becomes small as $t \to 0$; a relative error would therefore be inflated by division by a vanishing signal. In absolute terms the error of the learned transform remains nearly constant, while the coefficient error decreases together with the magnitude of the coefficient field. The order two coefficients should be zero for the affine additive problem; their small but nonzero values measure the residual nonlinear curvature already visible in Figure 2.

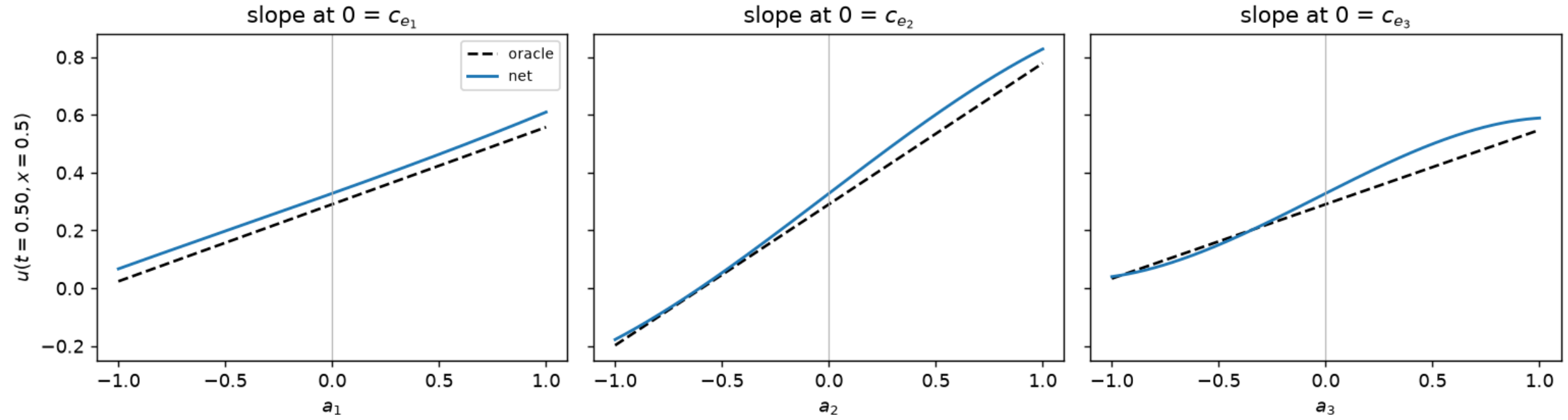

Figure 2: One dimensional slices of the learned restricted $S$-Transform. In each panel only one transform coordinate $a_k$ varies, while the other coordinates are set to zero. The dashed line is the closed form affine transform and the solid line is the learned map. The slope at $a = 0$ gives the corresponding first order chaos coefficient.

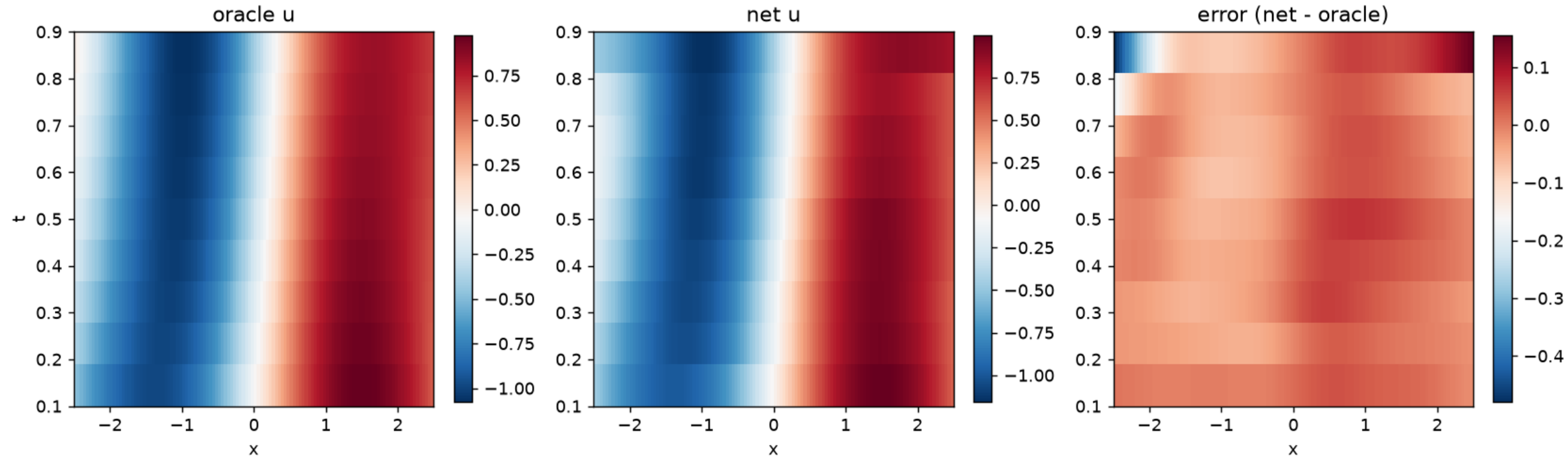

Figure 3: Learned restricted transform at the fixed transform parameter $a = (0.5, -0.5, 0.5)$. The left panel shows the closed form solution, the middle panel shows the neural approximation, and the right panel shows the signed error. The dominant error is located near the spatial boundary of the sampled region.

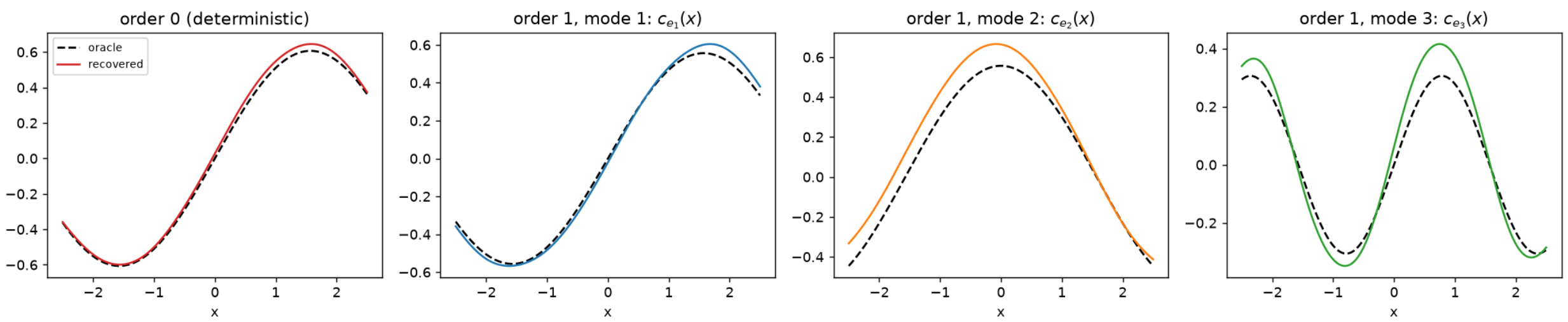

Figure 4: Recovered chaos coefficient fields at $t = 0.5$. Solid curves are the coefficients obtained by differentiating the learned transform at $a = 0$, while dashed curves are the closed form coefficients. The panels show the deterministic coefficient and the three retained first order noise modes.

| $t$ | $u$ RMSE | $\|u\|$ | $c_k$ RMSE | $\|c_k\|$ | order-2 RMSE |
|---|---|---|---|---|---|
| 0.9 | 0.057 | 0.52 | 0.087 | 0.50 | 0.020 |
| 0.7 | 0.056 | 0.52 | 0.080 | 0.43 | 0.020 |
| 0.5 | 0.055 | 0.52 | 0.071 | 0.34 | 0.018 |
| 0.3 | 0.043 | 0.57 | 0.051 | 0.23 | 0.017 |
| 0.1 | 0.025 | 0.63 | 0.023 | 0.09 | 0.022 |

Table 2: Absolute RMS error of the learned map and of the recovered first order coefficients across physical time $t$, with the corresponding signal magnitudes. The order two coefficients are zero in exact arithmetic.

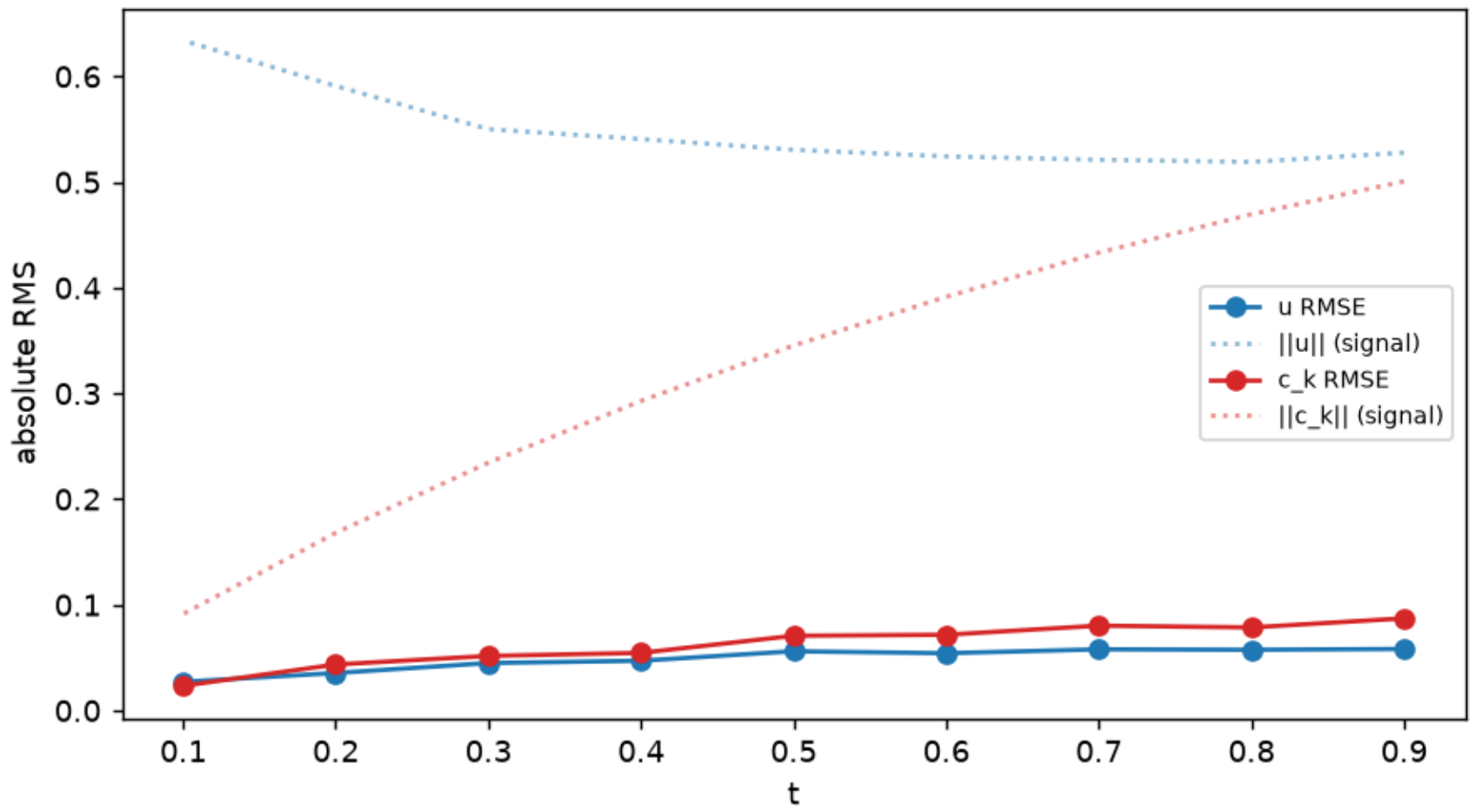

Figure 5: Absolute RMS errors and signal magnitudes as functions of time. The error of the learned transform is nearly flat, while the first order coefficient signal decreases toward $t = 0$. This explains why absolute errors are more informative than relative errors for this example.

Finally, Figure 6 reports the training loss for the backward learning scheme. The curves show rapid initial decrease and then fluctuate around a stable level. The step closest to the terminal condition starts from the largest loss, while the later backward steps benefit from warm starting and are already close to their final regime.

Recall that the time grid is coarse (ten steps), so part of the error may be attributed to time discretization and will shrink if the grid is refined. Because the heat equation with additive noise has a known explicit solution, our Example 2 validates the full NeSTR pipeline in a controlled setting: restrict, learn, differentiate, and reconstruct via (16). However, as the additive heat equation is affine linear, it does not test how the NeSTR methodology recovers high order coefficients in SPDEs which are nonlinear in noise. That is the role of the next example:

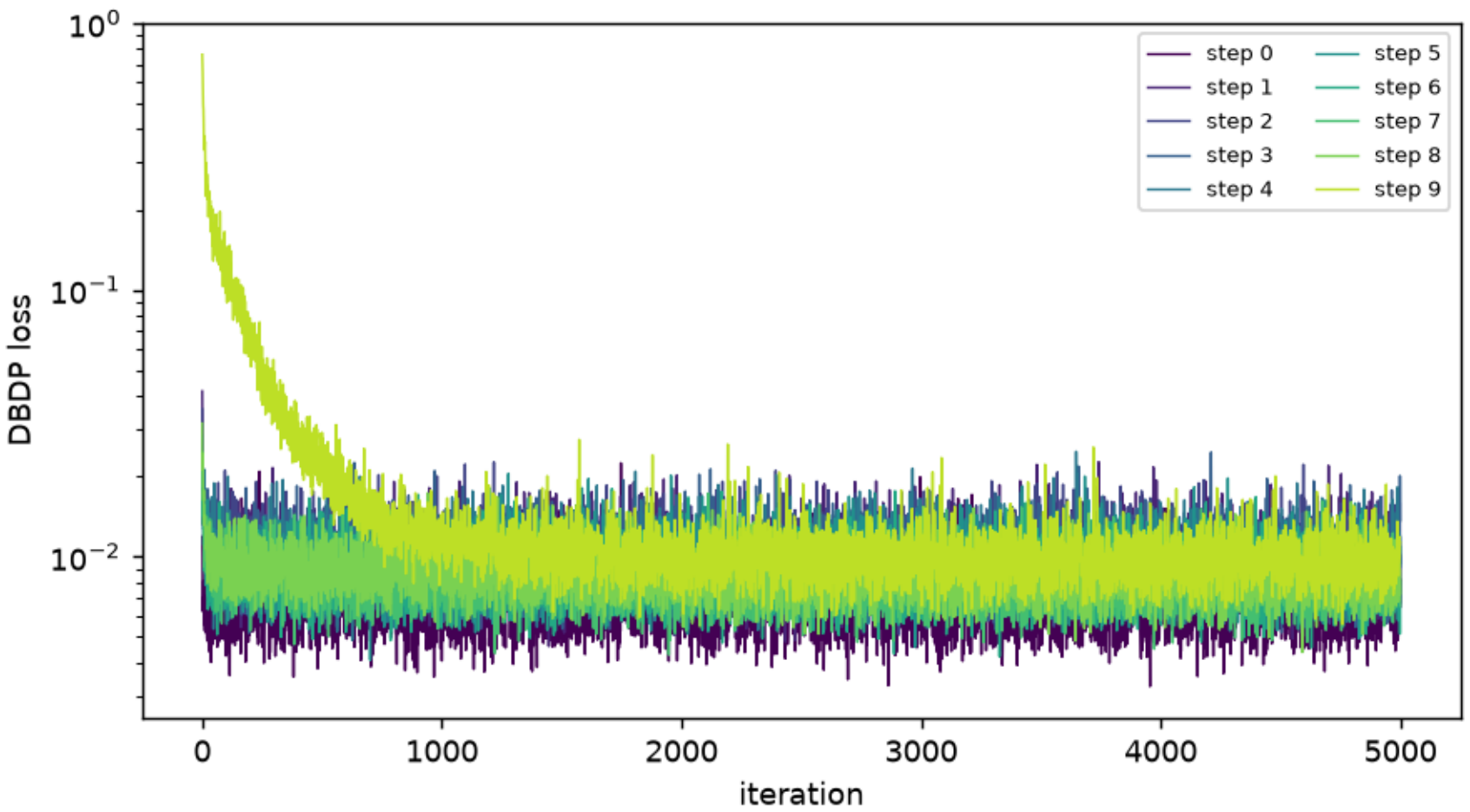

Figure 6: Training loss for the backward time stepping procedure. Each curve corresponds to one backward step.

**Example 3: Multiplicative Wick heat equation.** The third test is the multiplicative Wick heat equation

$$\begin{cases} \partial_t U(t,x) = \Delta U(t,x) + \sigma U(t,x) \diamond \dot{W}(t), & (t,x) \in (0,T] \times \mathbb{T}, \\ U(0,x) = u_0(x), & x \in \mathbb{T}, \end{cases} \tag{26}$$

with periodic boundary conditions. Here $W$ is temporal white noise, acting uniformly in space, and $\sigma > 0$ is a deterministic coupling strength. Note that the term $U(t,x) \diamond \dot{W}(t)$ is the same as the Skorohod interpretation of a stochastic integral, becoming a stochastic integral in the sense of Walsh [28] for adapted integrands (see also [14, 15]). Applying the $S$-Transform and using (1) gives the deterministic equation with potential

$$\partial_t u(t,x;h) = \Delta u(t,x;h) + \sigma u(t,x;h)h(t), \qquad u(0,x;h) = u_0(x).$$

After restriction to a finite set of noise modes, the potential becomes $\sigma h^{(N)}(a)$, and the corresponding Feynman–Kac representation has the form

$$u^{(N)}(t,x;a) = \mathbb{E}\left[u_0(X_t^x) \exp\left(\sigma \int_0^t h^{(N)}(a)(t-s)\, ds\right)\right],$$

where $X^x$ has generator $\Delta$ on the torus. Thus the driver is now multiplicative in the solution, namely $f = \sigma y h^{(N)}(a)$, while the rest of the numerical scheme is unchanged. In contrast to the additive case, $a \mapsto u^{(N)}$ is genuinely nonlinear, so chaos coefficients of every order are nonzero. Hence, this is an example that tests the recovery step (12) for higher-order derivatives.

To keep an exact target at all orders, we retain the normalized constant temporal mode in $H = L^2([0,1])$,

$$e_1(s) \equiv 1, \qquad h^{(1)}(a)(s) = a e_1(s) = a,$$

and choose the coupling $\sigma = \sqrt{2}$. The Gaussian coordinate $\xi_1 = W(e_1)$ is therefore standard normal, while the desired benchmark strength is carried by $\sigma$, not by a unnormalized noise mode. With this convention the integral in the Feynman–Kac exponent is

$$\sigma \int_0^t h^{(1)}(a)(t-s)\, ds = a\sqrt{2}\, t,$$

which is the source of the factor $\sqrt{2}\,t$ below. This normalization is chosen solely to produce a nontrivial exact benchmark with visible higher order chaos coefficients. Then the potential is deterministic and the restricted transform is

$$u^{(1)}(t,x;a) = e^{a\sqrt{2}\,t}(e^{t\Delta}u_0)(x).$$

Here $e^{t\Delta}$ is the heat semigroup on the torus; since $\Delta \sin x = -\sin x$, the benchmark data give $e^{t\Delta}u_0 = e^{-t}\sin x$, so all retained coefficients share this spatial profile. Consequently, the exact chaos coefficients associated with the retained constant temporal mode are

$$c_m(t,x) = \frac{(\sqrt{2}\,t)^m}{m!}(e^{t\Delta}u_0)(x), \qquad m = 0,1,2,\ldots.$$

We use $u_0 = \sin x$, $T = 1$, and $a \in [-1,1]$, and evaluate the recovered coefficients at the interior time $t = 0.5$. At this time the coupling $\sqrt{2}\,t \approx 0.71$ is large enough that the higher coefficients are visible: $c_2$ is about one quarter of $c_0$, and $c_3$ is about six percent of $c_0$.

Figure 7 displays the exact restricted transform at $a = 0.5$. In this case

$$u^{(1)}(t,x;0.5) = \exp\big((0.5\sqrt{2}-1)t\big)\sin x,$$

so the spatial sine profile is retained while its amplitude evolves exponentially in time. This surface provides a direct geometric contrast with the affine transform surface in the additive example.

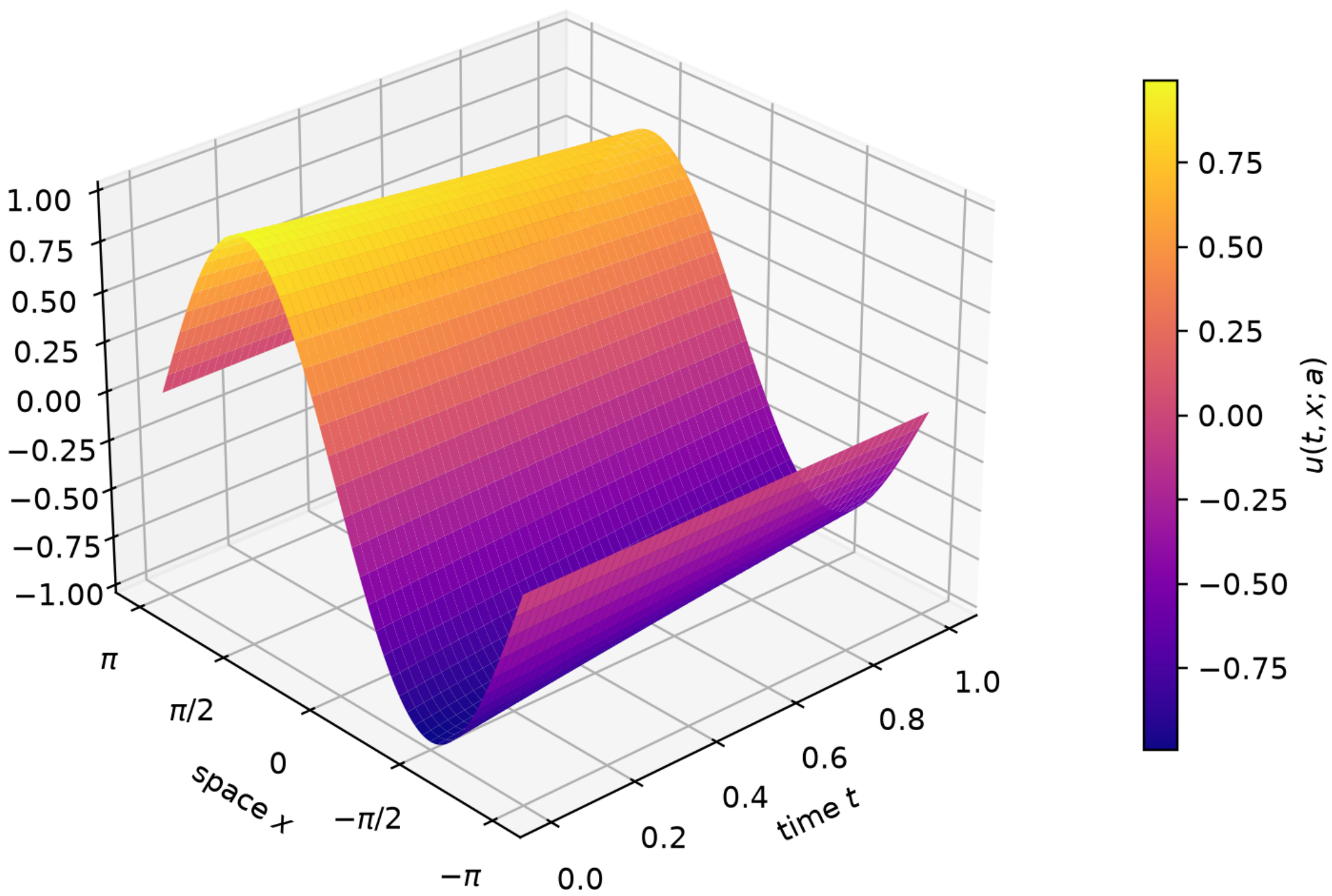

Figure 7: Exact multiplicative Wick restricted transform surface at $a = 0.5$ for $u_0(x) = \sin x$.

Figure 8 shows the recovered coefficient fields against the closed form. The first two coefficients are essentially exact, the second order coefficient is still captured, and the third order coefficient is already visibly less stable. This behavior is expected: the coefficients decay rapidly with order, so the signal to be recovered becomes smaller while the solve error stays roughly at the same scale.

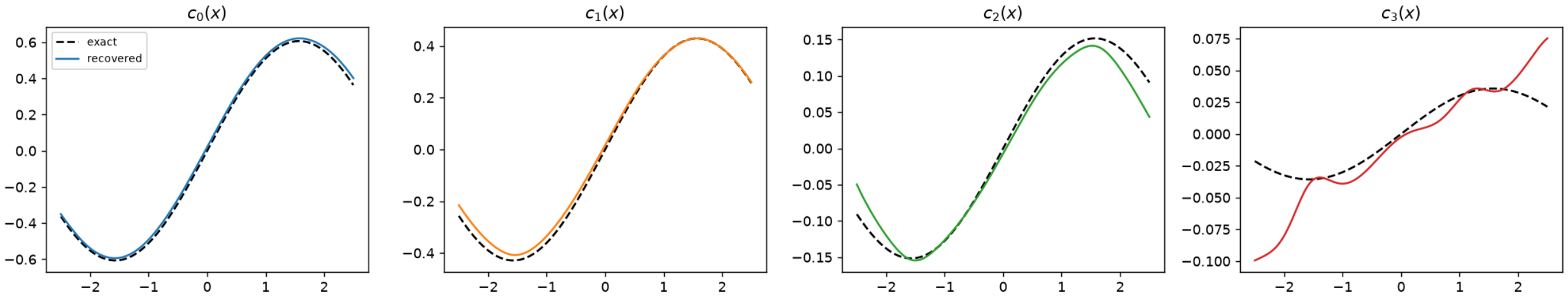

Figure 8: Recovered coefficient fields, solid curves, against the closed form, dashed curves, at $t = 0.5$, orders 0 through 3.

Recall the polynomial neural ansatz (23): the dependence on the transform coordinate $a$ is represented by a finite tensor polynomial expansion, here in the monomials $(a/A)^\nu$, while the coefficient functions $q^{(m)}_{\theta,\nu}(x)$ are represented by ordinary fully connected neural networks in the physical variable $x$. We compare this *polynomial spectral network* with a plain multilayer perceptron whose derivatives with respect to $a$ are computed by automatic differentiation. Table 3 reports the coefficient errors, and Figure 9 plots the same comparison by order. The polynomial spectral network is overall more accurate, and the advantage grows with the order. This is the expected benefit of representing the dependence on $a$ by an explicit polynomial and differentiating that polynomial exactly.

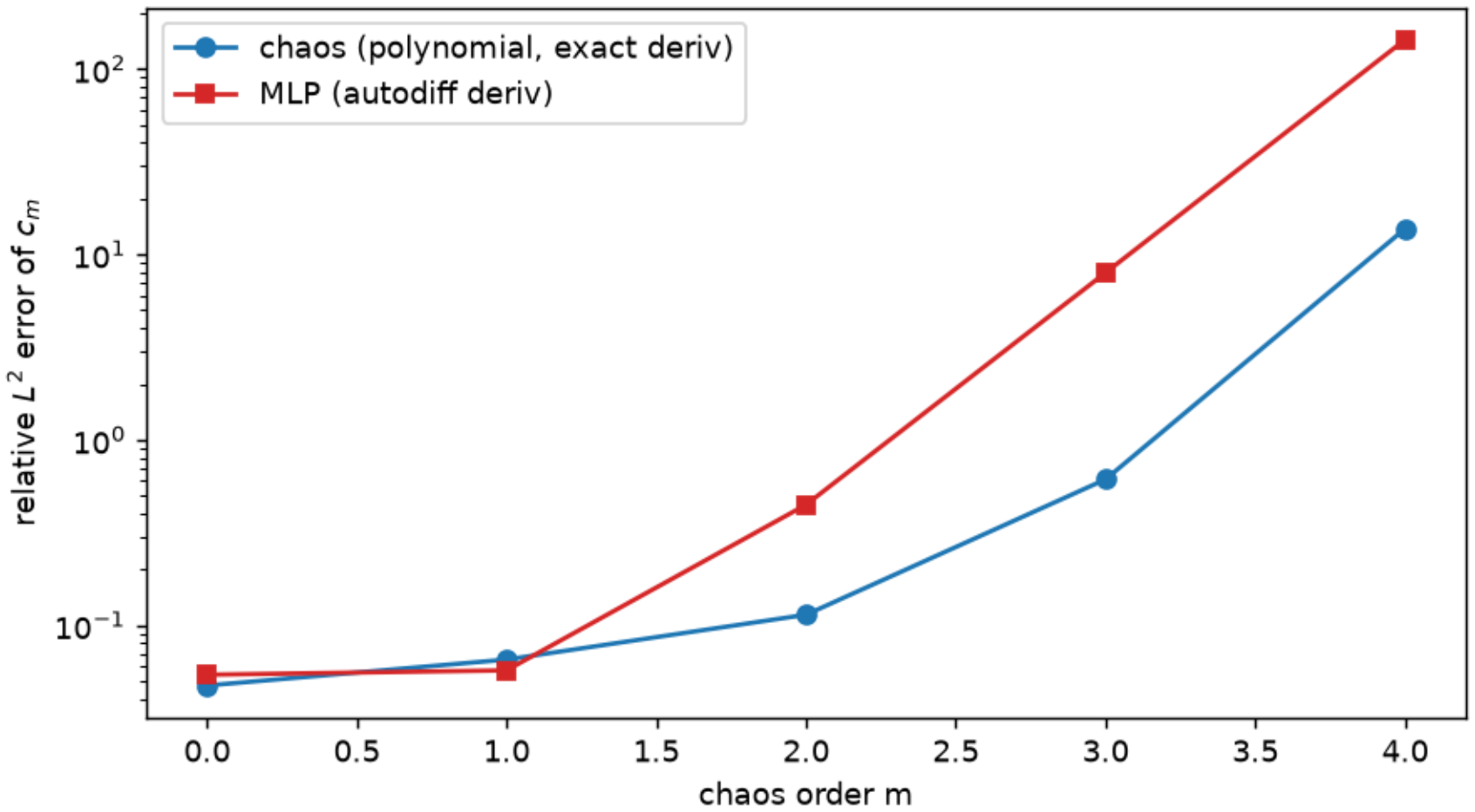

Figure 9: Coefficient recovery error by order, relative $L^2$ error on a log scale, for the polynomial spectral network and the multilayer perceptron.

The relative errors increase quickly with the order, but this should not be interpreted as a failure of the transform approximation alone. The exact coefficients shrink like $(\sqrt{2}\,t)^m/m!$, while the absolute errors in the last column of Table 3 remain comparable across several orders. Thus, a roughly fixed absolute error floor becomes a large relative error when the coefficient itself is very small.

Because the coefficients decay quickly, the reconstruction of $U$ through (16) can tolerate some noise in the highest orders. Table 4 and Figure 10 show the relative $L^2(\Omega)$ reconstruction error as a function of the chaos truncation order $M$.

| order $m$ | $\lvert c_m \rvert$ | spectral rel. $L^2$ | MLP rel. $L^2$ | spectral absolute error |
|---|---|---|---|---|
| 0 | 0.354 | 0.047 | 0.054 | 0.019 |
| 1 | 0.250 | 0.065 | 0.057 | 0.015 |
| 2 | 0.088 | 0.114 | 0.447 | 0.008 |
| 3 | 0.021 | 0.613 | 7.96 | 0.008 |
| 4 | 0.0037 | 13.7 | 143.6 | 0.044 |

Table 3: Coefficient recovery at $t = 0.5$, spectral network versus multilayer perceptron. The table reports relative $L^2$ errors, coefficient magnitudes, and the polynomial spectral network absolute error.

Two error sources contribute, and they separate exactly: the retained and the discarded chaos orders lie in orthogonal Hermite subspaces, so the squared relative error is the sum of a recovery part and a truncation part. The truncation part is known in closed form, the spatial factor $e^{t\Delta}u_0$ cancelling from the ratio, so that with the coupling $\rho := \sigma t$

$$\varepsilon_{\text{trunc}}(M) = \left( e^{-\rho^2} \sum_{m>M} \frac{\rho^{2m}}{m!} \right)^{1/2}. \tag{27}$$

This depends only on the retained index set, so it is common to both networks, and subtracting it isolates the two recovery mechanisms; Table 4 reports both. At $M = 0$ and $M = 1$ the total error is essentially this floor, so the networks necessarily agree there and those rows say nothing about either one. From $M = 2$ on the recovery part carries the whole difference. The same identity explains the optimal truncation $M = 3$: passing from $M = 2$ to $M = 3$ the floor falls by more than the recovery error grows, whereas passing to $M = 4$ it does not.

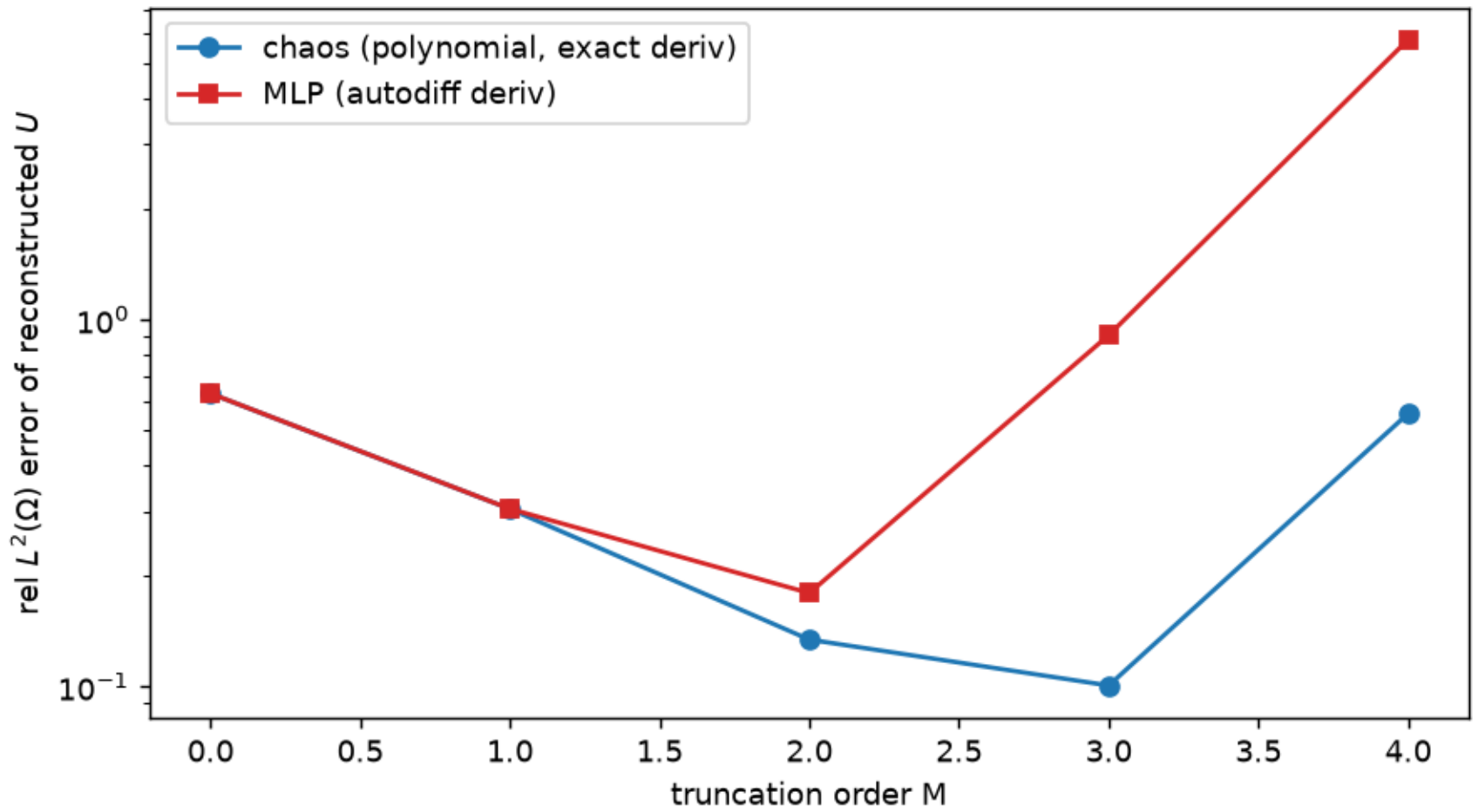

Figure 10: Relative $L^2(\Omega)$ reconstruction error of $U$ as a function of chaos truncation order $M$ for the polynomial spectral network and the multilayer perceptron.

Figure 11 shows the learned transformed field over the $(t, x)$ plane at a fixed value of $a$. The agreement with the closed form solution is good in the interior, with the largest discrepancies near the spatial edges and later times, where the field is largest. Figure 12 shows that training remains stable in the multiplicative case.

| $M$ | $\varepsilon_{\text{trunc}}(M)$ | total error | | recovery part only | |
|---|---|---|---|---|---|
| | | spectral | MLP | spectral | MLP |
| 0 | 0.627 | 0.628 | 0.629 | 0.030 | 0.047 |
| 1 | 0.300 | 0.305 | 0.305 | 0.053 | 0.053 |
| 2 | 0.120 | 0.134 | 0.180 | **0.060** | 0.134 |
| 3 | 0.042 | **0.101** | 0.906 | **0.092** | 0.905 |
| 4 | 0.013 | 0.554 | 5.78 | 0.554 | 5.78 |

Table 4: Relative $L^2(\Omega)$ error of the reconstructed stochastic field as a function of the chaos truncation order $M$, at $t = 0.5$. The truncation floor $\varepsilon_{\text{trunc}}$ of (27) is common to both networks; the last two columns are the remaining recovery error. Removing the floor raises the ratio between the two networks from 1.3 to 2.2 at $M = 2$ and from 9 to 10 at $M = 3$.

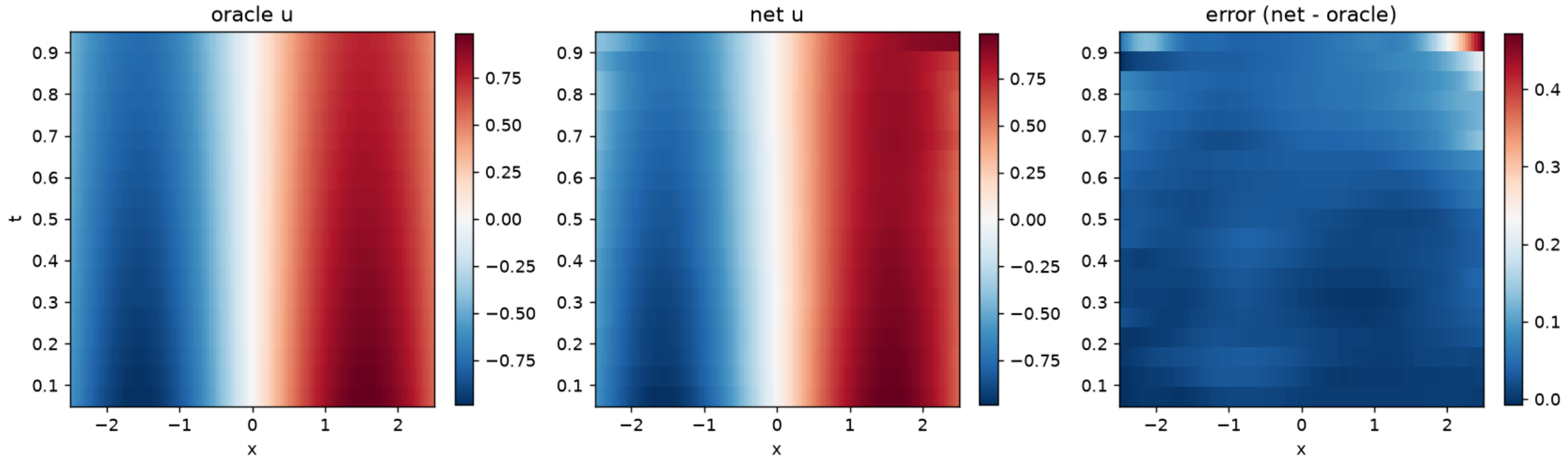

Figure 11: Multiplicative Wick field at $a = 0.5$: closed form solution, network, and signed difference.

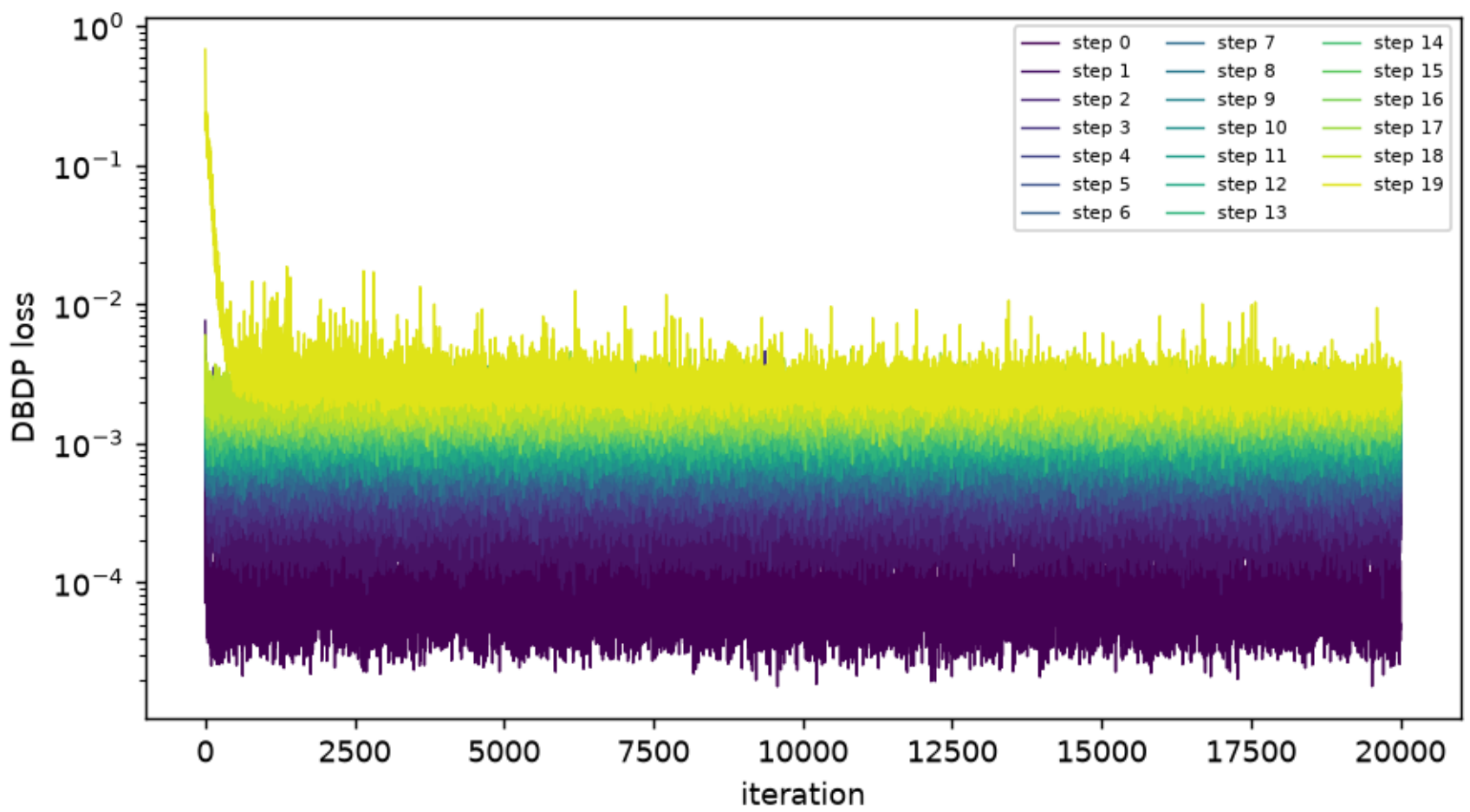

Figure 12: Training objective for the multiplicative Wick problem, shown by backward step. As in the additive experiment, only the step closest to the terminal datum starts far from its minimum.

For both networks the recovery error at order four exceeds the coefficient itself, so the usable truncation in this run is $M = 3$. This is the main limitation of the present configuration, and its mechanism is plain: the absolute recovery error is roughly constant across orders while the exact coefficients decay like $\rho^m/m!$, so beyond some order the signal falls below the error, and a larger coupling $\rho = \sigma t$ raises the high coefficients but by (27) also raises the truncation floor. As in the additive experiment, the time grid is still coarse. The balance between time discretization error and the accumulated per step network error remains an implementation parameter to tune.

We can also compare the numerical outcome with existing neural SPDE strategies at a qualitative level. Sample-wise deep solvers for SPDEs learn one realization or one noise path at a time [1], while modal space PINN methods learn stochastic modal dynamics [30] and recent Wiener chaos neural methods learn the chaos coefficients directly [23]. The present experiments show a different numerical behavior: NeSTR learns the restricted $S$-Transform once, and the chaos coefficients are then obtained by differentiation at the origin. In the additive test this recovers the expected affine transform and suppresses the spurious second order coefficients to a small error floor. In the multiplicative test the spectral representation resolves the reconstruction up to order $M = 3$, reaching a relative $L^2(\Omega)$ reconstruction error of about $10^{-1}$, whereas the plain multilayer perceptron becomes unreliable once unresolved higher order derivatives are included. Thus, compared with sample-wise or coefficient-wise approaches, the main numerical advantage observed here is not only the final field error, but the explicit and stable recovery of several chaos levels from a single deterministic transform approximation.

A further practical advantage is that the stochastic basis is reusable across equations. Once samples of the finite family $\{H_\alpha(\xi) : |\alpha| \leq M\}$ have been generated, they can be stored and used for every SPDE represented with the same Gaussian modes and chaos truncation. Each new SPDE contributes only its learned deterministic coefficient fields $c_\alpha(t, x)$. Combining different coefficient families with the same stored Hermite samples then produces sample paths or fields for several SPDEs without repeating the random basis simulation. This makes the separation between the deterministic coefficient computation and the universal random chaos basis not only conceptual, but also computationally useful.

# 5 Conclusion and outlook

We introduced NeSTR, a neural framework for solving Gaussian Wick type SPDEs by learning a finite dimensional restriction of the $S$-Transform. Rather than approximating stochastic trajectories or separate chaos coefficients, NeSTR learns the deterministic analytic map

$$(t, x, a) \longmapsto u^{(N)}(t, x; a),$$

whose derivatives at the origin recover the retained Wiener chaos coefficients. The stochastic solution is then obtained through an explicit truncated Hermite expansion. The distinctive contribution is therefore the transform as generating function formulation together with its derivative based inverse reconstruction.

The $S$-Transform first converts Wick products into ordinary products, producing a deterministic parametric PDE for learning. Proposition 3.9 then shows that, for fixed noise coordinates and truncation levels, the Hermite basis is universal across SPDEs sharing the same Gaussian noise; only the deterministic coefficient fields change. Theorem 3.13 complements this transfer principle by separating the final reconstruction error into noise mode, chaos truncation, deterministic solver, and neural transform contributions.

The additive and multiplicative heat equation experiments demonstrate accurate coefficient recovery and show that tensor polynomial representations of the transform variable provide stable derivatives at the origin.

The present results are limited to finite noise and chaos truncations and controlled equations with closed form references. Future work should establish quantitative derivative reconstruction rates, compare NeSTR directly with sample wise and coefficient wise neural methods, and address nonlinear Wick interactions, higher dimensional domains, and parameterized families of SPDEs such as those SPDEs considered in Holden *et al.* [15] beyond the class of stochastic heat equations we have treated here.